\documentclass[leqno,11pt]{amsart}

\usepackage{amsmath}
\usepackage{graphicx, color}
\usepackage{amscd}
\usepackage{amsfonts}
\usepackage{amssymb}
\usepackage{mathrsfs}
\usepackage{mathtools}

\newtheorem{thm}{Theorem}[section]

\newtheorem{lem}[thm]{Lemma}
\newtheorem{prop}[thm]{Proposition}

\newtheorem{rem}{Remark}[section]

\numberwithin{equation}{section}

\newcommand\be{\begin{equation}}
\newcommand\ee{\end{equation}}
\newcommand\R{\mathbb R}
\newcommand\Rnp{\mathbb R^n_+}
\allowdisplaybreaks
\def\eps{\varepsilon}

\title[Liouville theorem and superquadratic Hamilton-Jacobi equations]
{A Liouville-type theorem in a half-space and its\\ applications to the gradient blow-up behavior for superquadratic diffusive Hamilton-Jacobi equations}

\author{Roberta Filippucci}
\author{Patrizia Pucci}
\author{Philippe Souplet}

\begin{document}

\begin{abstract}
We consider the elliptic and parabolic superquadratic diffusive Hamilton-Jacobi equations:
$\Delta u+|\nabla u|^p=0$ and $u_t=\Delta u+|\nabla u|^p$,
with $p>2$ and homogeneous Dirichlet conditions.
For the elliptic problem in a half-space, we prove a Liouville-type classification, or symmetry result,
which asserts that any solution has to be one-dimensional.
This turns out to be an efficient tool to study the behavior of boundary gradient blow-up (GBU) solutions
of the parabolic problem in general bounded domains of~$\mathbb R^n$ with smooth boundaries.

Namely, we show that in a neighborhood of the boundary, at leading order, solutions display a global ODE type behavior
 of the form $u_{\nu\nu}\sim -u_\nu^p$, with domination of the normal derivatives
upon the tangential derivatives. 
This leads to the existence of a universal, sharp blow-up profile in the normal direction at any GBU point,
and  moreover implies that the behavior in the tangential direction is more singular.
A description of the space-time profile is also obtained.  
The ODE type behavior and its connection with the Liouville-type theorem can be considered as an analogue of the
well-known results of Merle and Zaag \cite{MZ98} for the subcritical semilinear heat equation,
with the significant difference that for the latter, the ODE behavior is in the time direction (instead of the normal spatial direction).

On the other hand, it is known that any GBU solution admits a weak continuation,
under the form of a global viscosity solution. As another consequence, we show that these viscosity solutions
{\it generically} lose boundary conditions after GBU. Namely, solutions without loss of boundary conditions after GBU
are exceptional and can be characterized as thresholds between global classical
solutions and GBU solutions which lose boundary conditions.
This result, as well as the  above GBU profile, were up to now essentially known only in one space-dimension.

Finally, in the case of elliptic Dirichlet problems, we deduce from our Liouville theorem an optimal Bernstein-type estimate,
which gives a partial improvement of a local estimate of P.-L.~Lions~\cite{Lions85}.

\smallskip
{\bf Keywords.} Diffusive Hamilton-Jacobi equations, Liouville-type theorem, gradient blow-up, final profile, loss of boundary conditions,
Bernstein-type estimates.
\end{abstract}

\maketitle

\section{Introduction and main results}
In this paper we consider superquadratic diffusive Hamilton-Jacobi equations, in both elliptic and parabolic settings.
Namely our goal is two-fold:

\begin{itemize}
\item[(i)] to establish a Liouville-type classification theorem for the elliptic problem in a half-space:
\be\label{main_pbLiouville}
\begin{cases}
-\Delta u =|\nabla u|^p, &x\in\Rnp,
\\ u(\tilde x,0)=0, & \tilde x\in\mathbb R^{n-1}, 
\end{cases}
\ee
where $\Rnp=\{(\tilde x,x_n)=(x_1,\dots,x_n)\in\mathbb R^n,\,\, x_n>0\}$ and $p>2$; 
\item[(ii)] to derive a number of applications of our Liouville-type theorem to the initial-boundary value problem:
\be\label{main_pb}\begin{cases}
u_t-\Delta u =|\nabla u|^p, & \mbox{in }\Omega\times(0,\infty),\\
u=0, &  \mbox{on }\partial\Omega\times(0,\infty),\\
u(\cdot,0) = u_0, & \mbox{in }\Omega,
\end{cases}
\ee
as well as for the inhomogeneous Dirichlet problem:
\be\label{ellinhomog}\begin{cases}
-\Delta u =|\nabla u|^p+f(x), &\mbox{in }\Omega,
\\ u=0, &\mbox{on }\partial\Omega.
\end{cases}
\ee
\end{itemize}
Throughout this article, $\Omega$ is a bounded domain of $\mathbb R^n$ with boundary
of class $C^{2+\mu}$ for some $\mu\in (0,1)$.
Let us denote by
$\nu_x$ the {\it inward} unit normal vector at any $x\in \partial\Omega$.

Let us begin with our Liouville-type classification, or symmetry, result. It asserts that any solution in a half-space is one-dimensional.

\begin{thm}\label{LiouvilleThm}
Let $p>2$ and let $u\in C^2(\Rnp)\cap C(\overline{\Rnp})$ be a solution of \eqref{main_pbLiouville}.
Then $u$ depends only on the variable $x_n$.
\end{thm}

As a consequence of Theorem~\ref{LiouvilleThm} and straightforward ODE analysis,
any solution of \eqref{main_pb} is thus given by
 either $u=0$ or $u=U_\alpha(x_n):=U_0(\alpha+x_n)- U_0(\alpha)$ 
for some $\alpha\ge 0$, where
\be\label{RefU0}
U_0(s)=c_ps^{1-\beta},\quad s\ge 0,
\qquad\hbox{ with } c_p={p-1\over p-2}(p-1)^{-1/(p-1)} =\frac{\beta^\beta}{1-\beta}.
\ee
Here and in the rest of the paper we define
$$\beta={1\over p-1}.$$
For future reference we also write
\be\label{RefU0prime}
U'_0(s)=d_ps^{-\beta},\quad s>0,
\qquad\hbox{ with } d_p=\beta^\beta=(1-\beta)c_p,
\ee
and we note that all solutions $U_\alpha$ for $\alpha>0$ are $C^1$, whereas 
$U'_0$ is singular at $s=0$ and $U_0$
displays the key H\"older exponent $1-\beta$.
We stress that Theorem~\ref{LiouvilleThm} does not assume $C^1$ regularity at the boundary for $u$,
and that this feature will be crucial in our applications.
We do not make any a priori assumption on the behavior of $u$ at infinity either.

\begin{rem}\label{rem0} {\rm 
(a) For the whole space case, it was proved in \cite{Lions85} that any classical solution
of $-\Delta u =|\nabla u|^p$ in $\R^n$ with $p>1$ has to be constant.
For the half-space problem \eqref{main_pbLiouville} in the subquadratic case $p\in (1,2]$,
a result similar to Theorem~\ref{LiouvilleThm} 
was proved in~\cite{PV}.
Our proof is based on a moving planes technique, combined with Bernstein type estimates
from~\cite{Lions85} and a compactness argument.
It is rather different from the proof in \cite{PV}, which relies on the existence of a finite limit as $x_n\to \infty$,
a property which does not hold in the superquadratic case. 

 Let us finally recall that Bernstein type gradient estimates go back to the early work \cite{Ber} and that the technique
was further developed in important papers such as \cite{La58}, \cite{Ar69}, \cite{Se69}, \cite{LY75}, \cite{Lions85}.
\smallskip 

(b) The Liouville-type theorem in \cite{PV} was motivated by the study of the so-called ``large solutions''
of elliptic equations with gradient terms, initiated in \cite{LL} in the framework of stochastic control
problems with state constraints (see also, e.g., \cite{BG}, \cite{LP1}, \cite{LP2}, \cite{AGQ}, \cite{FQS} and the references therein).
As for the question of 
one-dimensional symmetry of solutions of elliptic equations in a half-space,
it has also attracted much attention, especially for equations of the form $-\Delta u=f(u)$,
see, e.g., \cite{BCN}, \cite{FV}, \cite{Dup}, \cite{Fa15} and the references therein.}
\end{rem}

Let us turn to the applications of Theorem~\ref{LiouvilleThm} to the study of the parabolic problem~\eqref{main_pb}. 
Problem \eqref{main_pb} is locally well-posed for all $u_0\in X$,  with
$$X:=\{\phi\in C^1(\overline\Omega);\, \phi_{|\partial\Omega}=0\}.$$
Denoting by $T=T(u_0)$ the existence time of the unique maximal classical
solution $u$ of~\eqref{main_pb},
it is known that
$$\|u(\cdot,t)\|_\infty\le \|u_0\|_\infty,\quad 0<t<T,$$
as a consequence of the maximum principle, 
and that
$$T<\infty \Longrightarrow \lim_{t\to T^-}\|\nabla u(\cdot,t)\|_\infty=\infty.$$ 
This is called gradient blow-up (GBU) and it is also known that $T<\infty$ whenever $u_0$ is suitably large,
whereas solutions exist globally and decay to $0$ if $\|u_0\|_{C^1}$ is sufficiently small
(see e.g. \cite{ABG}, \cite{Alaa}, \cite{Sou}, \cite{HM}). 
The singular set, or GBU set, of $u$ is defined by
$$GBUS(u_0)=\Bigl\{x_0\in\overline\Omega;\, \limsup_{t\to T^-,\,x\to x_0} |\nabla u(x,t)|=\infty\Bigr\}$$ 
and the elements of $GBUS(u_0)$ are called GBU points.
It is known \cite{SZ} that
$$GBUS(u_0)\subset \partial\Omega.$$
More precisely, we have the following upper bound of Bernstein type:
\be\label{estBernstein1}
|\nabla u(x,t)|\le C(n,p)\delta^{-\beta}(x)+C(u_0)\quad\hbox{ in $\Omega\times [0,T)$,}
\ee
where
$$
\delta(x)={\rm dist}(x,\partial\Omega),\quad x\in\Omega,
$$
is the distance to the boundary (see \cite{SZ} and cf.~\cite{Lions85} in the elliptic case).
This also implies
\be\label{estBernstein2}
 -C(u_0) \delta(x)\le u(x,t)\le C(n,p)\delta^{1-\beta}(x)+C(u_0)\delta(x)\quad\hbox{ in $\Omega\times [ 0,T)$}
\ee
(the upper estimate follows by integrating \eqref{estBernstein1} in the normal direction; the lower estimate is immediate
since $u$ is a supersolution of the heat equation).
In view of \eqref{estBernstein1} and parabolic estimates, the solution $u$, which primarily belongs to $C^{2,1}(\overline\Omega\times (0,T))\cap C^{1,0}(\overline\Omega\times[0,T))$,
can be extended to a function $u\in C^{2,1}(\Omega\times (0,T])\cap C(\overline\Omega\times[0,T])$.
\smallskip

As a first consequence of Theorem~\ref{LiouvilleThm}, we have the following optimal Bernstein-type upper estimate,
which improves \eqref{estBernstein1}.
 Actually, the optimality will follow from Theorem~\ref{thmprofile1}.

\begin{thm} \label{Bernstein2}
Let $p>2$ and let $u_0\in X$ be such that $T(u_0)<\infty$.
For any
$\eps>0$ there exists $C=C(\eps,u_0)>0$ such that
\be\label{est1thmBernstein2}
|\nabla u(x,t)|\le(1+\eps)d_p \delta^{-\beta}(x)+C \quad\hbox{ in $\Omega\times [ 0,T]$,}
\ee
\be\label{est2thmBernstein2}
u(x,t)\le(1+\eps)c_p \delta^{1-\beta}(x)+C \delta(x) \quad\hbox{ in $\Omega\times [ 0,T]$,}
\ee
where the constants $c_p,d_p$ are given by \eqref{RefU0}, \eqref{RefU0prime}.
\end{thm}

In view of the next statements, let us introduce some notation. Set 
$$\Omega_\eps=\{x\in\Omega;\, \delta(x)<\eps\},\quad \eps>0.$$
Recall that, thanks to the regularity of $\Omega$, there exists $\delta_0>0$ such that
$$\hbox{for all $x\in\Omega_{\delta_0}$ there exists a unique
 point $P(x)\in\partial\Omega$ such that $|x-P(x)|=\delta(x)$.}$$
The point $P(x)$ is the projection of $x$ onto $\partial\Omega$. 
In this way, we can in particular extend the normal vector field
to the neighborhood $\Omega_{\delta_0}$ of the boundary, by setting, for all $x\in\Omega_{\delta_0}$:
\be\label{flowcoord}
\nu(x)=\nu_{P(x)},\qquad x=P(x)+\delta(x)\nu(x).
\ee
Also, for any unit vector fields $\xi,\zeta$, we denote the corresponding first and second order derivatives in space of $u$ by
$$u_\xi(x,t)=\nabla u(x,t)\cdot \xi(x),\qquad u_{\xi\zeta}(x,t)=\xi(x)D^2u(x,t)\zeta(x).$$

The following theorem describes a global behavior of the normal derivatives, 
and their dominance with respect to the tangential derivatives, for any GBU solution of~\eqref{main_pb}.

\begin{thm} \label{ODEbehavior}
Let $p>2$ and let $u_0\in X$ be such that $T(u_0)<\infty$.
For any $\eps>0$ there exists $C_\eps>0$ (possibly depending on $u_0$)  
such that
\be\label{boundnormal}
\bigl|u_{\nu\nu}+|u_\nu|^p\bigr|\le \eps|u_\nu|^p+C_\eps
\quad\hbox{ in $\Omega_{\delta_0}\times [T/2,T]$}
\ee
and
\be\label{boundtangential}
|u_{\tau\tau}|+|u_{\nu\tau}|+|u_\tau|^p\le \eps|u_\nu|^p+C_\eps
\quad\hbox{ in $\Omega_{\delta_0}\times [T/2,T]$,}
\ee
where $\tau$ is any tangential vector field (i.e. $\tau \perp \nu$ and $|\tau|=1$ in $\Omega_{\delta_0}$).
Moreover, we have  
\be\label{boundinf}
\inf_{\Omega_{\delta_0}\times (0,T)}u_\nu >-\infty.
\ee
\end{thm}

Theorem~\ref{ODEbehavior} has interesting consequences on the behavior near a GBU point.

\begin{thm} \label{thmprofile1}
Let $p>2$ and let $u_0\in X$ be such that $T(u_0)<\infty$. For each GBU point $a\in \partial\Omega$,
we have the following properties:
\smallskip

(i) {\it (Universal final blow-up profile in the normal direction)}
$$\lim_{s\to 0} s^\beta \nabla u(a+s\nu_a,T)=d_p\nu_a.$$

(ii) {\it (More singular behavior in the tangential direction)}
\be\label{moresingul}
\lim_{x\to a,\, x\in\partial\Omega} |x-a|^\beta u_\nu(x,T)=\infty.
\ee

(iii)  {\it (Continuity of $u_\nu$ 
with values in $[0,\infty]$)} As $t\to T$ and $x\to a$, the normal derivative~$u_\nu$  
(and hence $|\nabla u|$) blows up in the strong sense:
\be\label{contnabla}
\lim_{t\to T^-,\,x\to a} u_\nu=\infty.
\ee

(iv) {\it (Space-time behavior)} We have
\be\label{spacetime1}
\lim_{t\to T^-,\,x\to a} \frac{u_{\nu\nu}}{u_\nu^p} =-1
\ee
and consequently, for each $\eps\in (0,1)$ there exists $\eta\in(0,\delta_0)$ such that
\be\label{spacetime2}\hspace{-.5truecm}
 \ \Bigl[u_\nu^{1-p}(P(x),t)+(1+\eps)(p-1)\delta(x)\Bigr]^{-\beta}\le u_\nu(x,t)\le \Bigl[u_\nu^{1-p}(P(x),t)+(1-\eps)(p-1)\delta(x)\Bigr]^{-\beta}
 \ee
 for all $(x,t)\in B_\eta(a)\times (T-\eta,T)$.
 Moreover, we have
\be\label{spacetime3}
\bigl|u_\nu^{1-p}(z,t)-u_\nu^{1-p}(a,t)\bigr|\le o(|z-a|)\quad\hbox{ as $t\to T$ and $z\to a$, with $z\in\partial\Omega$.}
\ee
\end{thm}

\goodbreak
Theorems~\ref{ODEbehavior} and \ref{thmprofile1} show that any GBU solution follows a global ODE-like behavior in the normal direction.
More precisely, in the singular region, the dominating terms in the PDE $u_t-\Delta u=|\nabla u|^p$
are the normal derivatives $u_\nu$ and $u_{\nu\nu}$, with
$$u_{\nu\nu}\sim -|u_\nu|^p,$$
whereas all other derivatives are of lower order,
as illustrated by the following scheme (in two space dimensions):
\def\boxitA#1{\leavevmode\hbox{\vrule\vtop{\vbox{\kern.33333pt\hrule
   \kern5pt\hbox{\kern3pt\vbox{#1}\kern3pt}}\kern5pt\hrule}\vrule}}
 \def\boxitB#1{\leavevmode\hbox{\vrule\vtop{\vbox{\kern.33333pt\hrule
   \kern2pt\hbox{\kern3pt\vbox{#1}\kern3pt}}\kern3pt\hrule}\vrule}}
\setbox1=\vbox{
\hsize=6mm
\baselineskip=14pt \parindent=0mm
$u_{\nu\nu}$}
\setbox2=\vbox{
\hsize=4.5mm
\baselineskip=14pt \parindent=0mm
$u_\nu^2$}
$$u_t-\boxitA{\hbox{\box1}}-u_{\tau\tau} =\bigl(\boxitB{\hbox{\box2}}+u_\tau^2\,\bigr)^{p/2}.$$
For the tangential parts of the gradient and of the Laplacian, this is stated in \eqref{boundtangential}.
As for the time-derivative, we actually have the well-known uniform bound
\be\label{boundutsol}
M:=\sup_{\Omega\times [T/2,T]} |u_t| <\infty
\ee
(this follows from the maximum principle applied to $u_t$, see e.g. \cite{SZ}).
We stress that the value $u_\nu(x,T)=\infty$ in~\eqref{moresingul} is allowed,
so that the statement applies to both isolated or non-isolated GBU points.

\smallskip 
We next turn to the post blow-up behavior. It is known \cite{BaLio04} that problem \eqref{main_pb} admits a unique, continuous
global viscosity solution $\bar u\in C(\overline\Omega\times [0,\infty))$, which extends the maximal classical solution after $t=T$.
It is actually a classical solution in $\Omega\times (0,\infty)$,
namely, $\bar u\in C^{2,1}(\Omega\times (0,\infty))$,
but the homogeneous boundary conditions have to be understood in the generalized viscosity sense (or state constraints)
and need not be satisfied in the usual sense.
The solution $\bar u$ can also be obtained by monotone approximation of problems with truncated nonlinearities
(see \cite{PZ_aihp12}, \cite{PS_jmpa18} and the references therein).
Moreover, as shown in \cite{PS_aihp17}, \cite{QR18}, the global weak solution~$\bar u$ may lose boundary conditions after gradient blow-up, i.e.
$$\sup_{t>T,\, x\in\partial\Omega}\bar u(x,t)>0.$$
However, it was shown in \cite{PS_aihp17} that there are also solutions which never lose the boundary conditions after gradient blow-up.
In particular, for any nontrivial 
$\phi\in X_+$, where $X_+=\{\phi\in X;\, \phi\ge 0\}$, it is shown in \cite{PS_aihp17} that
$$\lambda^*:=\sup\{\lambda >0;\, T(\lambda\phi)=\infty\} \in (0,\infty),$$
that $T(\lambda^*\phi)<\infty$ and that the corresponding solution does not lose boundary conditions.
In one space dimension, it was moreover proved in \cite{PS_aihp17} that
solutions without loss of conditions are exceptional: they constitute thresholds between global classical solutions and GBU solutions
with loss of boundary conditions. As a consequence of Theorems~\ref{Bernstein2} and \ref{thmprofile1}, 
we can show that this threshold property remains true in any space dimension.

\begin{thm}  \label{thmLBC}
Let $p>2$ and let $u_0,v_0\in X_+$ be such that $T(u_0)<\infty$ and $v_0\not\equiv u_0$.
Denote by $\bar u,\bar v$ the corresponding unique global viscosity solutions of~\eqref{main_pb}.

\begin{itemize}
\item[(i)] If $v_0\le u_0$ and no loss of boundary conditions occurs for $\bar u$, then $T(v_0)=\infty$.

\smallskip

\item[(ii)] If $v_0\ge u_0$, then $T(v_0)<T(u_0)<\infty$ and $\bar v$ loses boundary conditions
before $T(u_0)$, i.e., there exists $t\in(T(v_0),T(u_0))$ such that
$$\max_{x\in\partial\Omega} \bar v(x,t)>0.$$
\end{itemize}
\end{thm}

See \cite{PZ_aihp12}, \cite{PS_aihp17},~\cite{PS_jmpa18}
for further results on the behavior of the viscosity solution $\bar u$ for $t\ge T$.
We refer to \cite{FL94}, \cite{fl19} and the references therein for results on the continuation
after GBU for some other one-dimensional parabolic problems.
\medskip

Let us finally briefly consider the inhomogeneous elliptic problem~\eqref{ellinhomog}.
It is shown in \cite[Th\'eor\`eme~IV.1]{Lions85} that any (local) solution $u\in C^2(\Omega)$
of
\be\label{locEllipt}
-\Delta u =|\nabla u|^p+f(x)\quad\hbox{ in }\Omega,
\ee
with $f\in W^{1,\infty}(\Omega)$, satisfies the Bernstein-type estimate
\be\label{BernsteinPLL}
|\nabla u(x)|\le C(n,p) \bigl[\delta^{-\beta}(x)+\|f\|_{W^{1,\infty}}\bigr] \quad\hbox{ in $\Omega$.}
\ee
As a consequence of Theorem~\ref{LiouvilleThm},
we obtain the following optimal estimate,
 which is a partial improvement of the result in~\cite{Lions85}
for the case when the boundary value problem~\eqref{ellinhomog}
is considered instead of the local equation~\eqref{locEllipt}.

\begin{thm}  \label{thmEllipt}
Let $p>2$ and $f\in W^{1,\infty}(\Omega)$ with $\|f\|_{W^{1,\infty}}\le M$.
Assume that $u\in C^2(\Omega)\cap C(\overline\Omega)$
 is a solution of~\eqref{ellinhomog}. Then for every $\eps>0$
there exists $C=C(\eps,M)>0$ such that
$$ 
|\nabla u(x)|\le(1+\eps)d_p \delta^{-\beta}(x)+C \quad\hbox{ in $\Omega$.}
$$
\end{thm}

\begin{rem}\label{rem1} {\rm
As far as we know, this paper provides the first study of the spatial GBU behavior and final profiles valid
for general solutions of~\eqref{main_pb} in all space dimensions.
Previously, the behavior was known only in the one-dimensional case, see~\cite{CG96}, \cite{ARS04}, \cite{GuoHu},
or for domains with some symmetry \cite{PS_imrn16}.
\smallskip

Theorems \ref{LiouvilleThm} and \ref{ODEbehavior} can be seen as the analogues of
the well-known results of Merle and Zaag \cite{MZ98} (see also \cite{MZ00})
concerning the subcritical nonlinear heat equation
\be\label{NLH}
u_t-\Delta u=|u|^{p-1}u
\ee
with $1<p<(n+2)/(n-2)_+$.
As a key difference, the ODE behavior in \eqref{NLH} is in the time direction for $u$,
whereas the ODE behavior in \eqref{main_pb} is in the
spatial normal direction for $\nabla u$.
Namely, the Liouville-type theorem in~\cite{MZ98} states that any ancient solution of \eqref{NLH} with self-similar temporal decay at $-\infty$
must depend on the time-variable only. This is then used to show that
blow-up solutions of \eqref{NLH} satisfy
$$|u_t-|u|^{p-1}u|\le \eps|u|^p+C_\eps$$
(see also \cite{MizNeu}, \cite{fi_aihp14}, \cite{GuoSou2} for related results based on the Liouville theorem in \cite{MZ98}).

\smallskip

The proof of Theorem~\ref{ODEbehavior} relies on Theorem~\ref{LiouvilleThm},
combined with suitable rescaling and compactness arguments.
It follows the general strategy of~\cite{MZ98} (see also \cite{GuoSou2}), but with notable differences.
First, whereas the solutions of \eqref{NLH} considered in \cite{MZ98} blow up only at interior points of the domain $\Omega$
(as a consequence of a convexity assumption on $\Omega$),
GBU for \eqref{main_pb} occurs at the boundary. Due to this, we have to deal with rather delicate
boundary estimates in our rescaling procedures and in the preliminary nondegeneracy properties,
relying in particular on flow coordinates (cf.~\eqref{flowcoord}).
Moreover, the nondegeneracy properties require different arguments from those in~\cite{MZ98},
due to the lack of variational structure of problem~\eqref{main_pb}. 
Also, instead of using type~I temporal estimates from \cite{GK87} for \eqref{NLH} as basic a priori estimates,
we rely on the spatial Bernstein type estimate~\eqref{estBernstein1}.
\smallskip

As another qualitative difference with \cite{MZ98}, we note that our results on the parabolic problem \eqref{main_pb}
are derived from an {\it elliptic} Liouville-type theorem.
This is allowed by the above mentioned bound \eqref{boundutsol} on  $u_t$ in \eqref{main_pb},
so that the time derivative vanishes in rescaling limits.
Let us stress that the apparently simplifying mechanisms
\eqref{boundnormal}-\eqref{boundtangential} and \eqref{boundutsol}
 are far from making the dynamics of the equation trivial.
Indeed, they are not sufficient to provide  complete information on the transition (in space and/or in time) between
the singular and regular parts of the solution
and on the corresponding transition speeds (time rate of GBU and tangential space profile near an isolated GBU point of the boundary).
These questions are delicate; see Remark~\ref{rem1b}(b) 
for results in that direction.}
\end{rem}

\begin{rem}\label{rem1b} {\rm 
(a) Ancient solutions of \eqref{main_pb} in $\R^n\times(-\infty,0)$ have been
studied in~\cite{SZ}.
The half-space case $\Rnp\times(-\infty,0)$, which is a topic of possible independent interest, will be studied in a forthcoming paper.
\smallskip

(b) As mentioned above,
Theorem~\ref{ODEbehavior} is not sufficient to determine the sharp final GBU profile in the tangential direction
near isolated GBU points and we suspect that the latter is not universal but that various tangential profiles may exist, depending on the solution.
However, Theorem~\ref{thmprofile1}(ii) shows that the profile is always anisotropic: it is more singular in the tangential direction than in the normal one.

\smallskip
In some very special cases, more precise information on the final GBU profile in the tangential direction
can be found in \cite{PS_imrn16}. Namely, for $n=2$ and $p\in (2,3]$, under suitable symmetry
assumptions on the domain $\Omega$ and initial data $u_0$,
and assuming that $\Omega$ coincides with the half-plane $\{(x,y);\, y>0\}$ near the origin,
we have single-point GBU at the origin, with the final profile
$$u_\nu(x,y,T)\sim d_p\Bigl[y+C |x|^{2(p-1)/(p-2)}\Bigr]^{-\beta},\quad\hbox{ as $x,y\to 0$.}$$
For results on the GBU set, especially sufficient conditions ensuring single-point GBU, see \cite{LS10}, \cite{Esteve19}.
As for the time rate $\|\nabla u(\cdot,t)\|_\infty$ of GBU, it remains an open problem in general,
and so is the time behavior of $u_\nu(a,t)$ -- cf.~\eqref{spacetime3}.
Nevertheless the rate is known to be always non self-similar, unlike for \eqref{NLH}, and of type~II,
with a lower bound $\|\nabla u(\cdot,t)\|_\infty\ge C(T-t)^{-1/(p-2)}$;
see \cite{CG96}, \cite{GuoHu}, \cite{ZL13}, \cite{PS_jmpa18} for results on the GBU rate.
\smallskip

(c) Similar to \eqref{boundtangential}, we also have the same bound for the mixed tangential derivative:
\be\label{boundtangential2}
|u_{\tau\hat\tau}|\le \eps|u_\nu|^p+C_\eps
\quad\hbox{ in $\Omega_{\delta_0}\times [T/2,T]$},
\ee
where $\tau,\hat\tau$ are any vector fields such that $(\nu, \tau,\hat\tau)$ is orthonormal in $\Omega_{\delta_0}$
(see the proof of Theorem~\ref{ODEbehavior}).
}
\end{rem}

\bigskip

The outline of the rest of the paper is as follows.
In Section~2 we establish Theorem~\ref{LiouvilleThm}.
In Sections~3 and 4, we use Theorem~\ref{LiouvilleThm} to prove Theorems~\ref{Bernstein2}
and \ref{thmEllipt}, respectively. 
Theorem~\ref{ODEbehavior} is next proved in Section~5.
In Section~6 we deduce Theorem~\ref{thmprofile1} from Theorem~\ref{ODEbehavior}, and Theorem~\ref{thmLBC}
from Theorems~\ref{Bernstein2} and \ref{thmprofile1}.
Finally, in the appendix, we provide the proof of Proposition~\ref{lemCentral},
a technical result which is used in the proof of Theorem~\ref{ODEbehavior}.

\section{Proof of Theorem~\ref{LiouvilleThm}}

It is based on a moving planes argument combined with Bernstein type estimates.

\begin{proof}
 Write $x=(\tilde x,y)\in \mathbb R^{n-1}\times  [0,\infty)$ and
fix any $h\in \mathbb R^{n-1}\setminus\{0\}$.  Let
$$v(\tilde x,y)=u(\tilde x+h,y)-u(\tilde x,y),\qquad (\tilde x,y)\in \mathbb R^{n-1}\times [0,\infty).$$
It suffices to show that $v\equiv 0$.
Assume for contradiction that
\be\label{eqContrad}
\sigma:=\sup_{\mathbb R^n_+} v>0
\ee
(the case $\inf_{\mathbb R^n_+} v<0$ is similar). 
By the Bernstein estimate in \cite{Lions85}, we have
\be\label{Bernsteiny}
|\nabla u(\tilde x,y)|\le C(n,p)y^{-\beta},\quad\hbox{for all $(\tilde x,y)\in \mathbb R^{n-1}\times (0,\infty)$}
\ee
(more precisely, we apply Case 2  
of Theorem~IV.1 and Remark p.~250 in \cite{Lions85},
with $C_1=C_2=0$, to the function $-u$).
It follows that
$$|v(\tilde x,y)|\le C(n,p)|h|y^{-\beta},\quad\hbox{for all $(\tilde x,y)\in \mathbb R^{n-1}\times (0,\infty)$.}$$
 Hence $|v|\le \sigma/2$ for $y\ge A$ large.
Therefore
\be\label{eqContrad2}
\sigma=\sup_{\mathbb R^{n-1}\times (0,A)} v.
\ee
On the other hand, $v$ satisfies the equation
\be\label{eqConv}
-\Delta v=a(\tilde x,y)\cdot\nabla v \ \text{ in } \mathbb R^n_+,
\ee
where
$$a(\tilde x,y):=\int_0^1 G\bigl(s\nabla u(\tilde x+h,y)+(1-s)\nabla u(\tilde x,y)\bigr)\,ds,\qquad G(\xi)=p|\xi|^{p-2}\xi,$$
Observe that, as a consequence of \eqref{Bernsteiny}, the function $a$ is bounded for $y$ bounded away from~$0$,
hence in particular on compact subsets of $\mathbb R^n_+$.
By the strong maximum principle applied to~\eqref{eqConv}, it follows that
 the solution $v$ cannot achieve any local maximum in $\mathbb R^n_+$.
 Otherwise $v$ would be constant, and we assumed the contrary in view of \eqref{eqContrad} and the fact that $v(\tilde x,0)=0$.
In particular, $\sigma$ in \eqref{eqContrad2} is not attained and
there exists a sequence $(\tilde x_j,y_j)\in \mathbb R^{n-1}\times (0,A)$ with $|\tilde x_j|\to\infty$
such that $v(\tilde x_j,y_j)\to\sigma$.

\smallskip
Next define
$$u_j(\tilde x, y)=u(\tilde x_j+\tilde x,y),\qquad (\tilde x,y)\in \mathbb R^{n-1}\times [0,\infty),$$
and note that
\be\label{supduj}
\sup_{(\tilde x, y)\in\mathbb R^n_+} \bigl(u_j(\tilde x+h, y)-u_j(\tilde x, y)\bigr)=\sup_{\mathbb R^n_+} v=\sigma
\ee
and
\be\label{limsigma}
u_j(h,y_j)-u_j(0,y_j) =v(\tilde x_j,y_j) \to \sigma,\quad\mbox{ as }j\to\infty.
\ee
Since $u_j$ is a solution of \eqref{main_pb} in $\mathbb R^n_+$,
 it satisfies the (uniform) Bernstein estimate \eqref{Bernsteiny}, namely:
$$|\nabla u_j(\tilde x,y)|\le C(n,p)y^{-\beta}, \quad\hbox{for all $(\tilde x,y)\in \mathbb R^{n-1}\times (0,\infty)$.}$$
Owing to $u_j(\tilde x,0)=0$, by integration in the $y$ direction, we also have
\be\label{Bernsteiny_uj}
|u_j(\tilde x,y)|\le C_1(n,p)y^{1-\beta},\quad\hbox{for all $(\tilde x,y)\in \mathbb R^{n-1}\times [0,\infty)$}
\ee
 and for all $j$.
It then follows from interior elliptic estimates that $(u_j)_j$ is relatively compact in $C^2_{\mathrm{loc}}(\mathbb R^n_+)$.
Therefore, some subsequence of $(u_j)_j$ converges in that topology to a solution $U\in C^2(\mathbb R^n_+)$ of $-\Delta U =|\nabla U|^p$.
As a consequence of \eqref{Bernsteiny_uj}, we also have
$U\in C(\overline{\mathbb R^n_+})$ and $U(\tilde x,0)=0$.
Moreover, we may assume that $y_j\to y_\infty\in [0,A)$ and we get
\be\label{supU}
U(h,y_\infty)-U(0,y_\infty)= \sigma,
\ee
owing to \eqref{limsigma}, which implies $y_\infty>0$.

\smallskip
 Put now
$$V(\tilde x,y)=U(\tilde x+h,y)-U(\tilde x,y),\qquad (\tilde x,y)\in \mathbb R^{n-1}\times [0,\infty).$$
It follows from \eqref{supduj} and \eqref{supU} that $\sigma=\sup_{\mathbb R^n_+} V=V(0,y_\infty)$.
But $V$ satisfies
$$-\Delta V=A(\tilde x,y)\cdot\nabla V,$$
where
$$A(\tilde x,y):=\int_0^1 G\bigl(s\nabla U(\tilde x+h,y)+(1-s)\nabla U(\tilde x,y)\bigr)\,ds$$
is bounded on compact subsets of $\mathbb R^n_+$.
This contradicts the strong maximum principle and completes the proof.
\end{proof}

\section{Proof of Theorem~\ref{Bernstein2}}

The proof of Theorem~\ref{Bernstein2} is based on the Liouville-type Theorem~\ref{LiouvilleThm} and on a suitable rescaling argument.
By the same ideas, one also obtains the following proposition, which states that tangential derivatives
are of lower order than normal derivatives in terms of the distance to the boundary.
This will be an important preliminary step for the proof of Theorem~\ref{ODEbehavior} in Section~5.

\begin{prop}[] \label{lem1}
Let $p>2$ and let $u_0\in X$ be such that $T(u_0)<\infty$.
For each $\eps>0$ there exists a constant $C_\eps>0$ such that,
for any vector fields $\tau,\hat\tau$ such that $\tau,\hat\tau \perp \nu$
and $|\tau|=|\hat\tau|=1$ in $\Omega_{\delta_0}$, we have
\begin{equation}
\label{ConclLemGrad}
|\nabla u\cdot\tau|^p+|\tau (D^2u)\nu|+|\tau (D^2u)\hat\tau|\le \eps  \delta^{-\beta p}+C_\eps
\quad\hbox{in $\Omega_{\delta_0}\times[T/2,T)$.}
\end{equation}
\end{prop}

\medskip

{\it Proof of Theorem~\ref{Bernstein2} and of Proposition \ref{lem1}.}
First recall that $u$ satisfies the following Bernstein-type estimates:
\begin{eqnarray}
|\nabla u(x,t)|&\le& C\delta^{-\beta}(x)\quad\hbox{ in $\Omega\times [0,T)$,} \label{estBernstein1b} \\
|u(x,t)|&\le& C\delta^{1-\beta}(x)\quad\hbox{ in $\Omega\times [0,T)$,} \label{estBernstein2b}  
\end{eqnarray}
for some constant $C>0$ depending only on $\Omega, p, u_0$; cf.~\eqref{estBernstein1}-\eqref{estBernstein2}.

Assume that either \eqref{est1thmBernstein2} or \eqref{ConclLemGrad} fails. Then there exist $c>0$,
a sequence of points $(x_j,t_j)$ in~$\Omega_{\delta_0}\times[T/2,T)$
and unit vectors $\tau_j, \hat\tau_j\perp \nu(x_j)$ such that either
\begin{equation}
\label{HypContradDelta0a}
|\nabla u(x_j,t_j)|\ge (1+c)d_p\, \delta^{-\beta}(x_j)+j
\end{equation}
or
\begin{equation}
\label{HypContradDelta0}
\Bigl\{|\nabla u\cdot\tau_j|^p+|\tau_j (D^2u)\nu|+|\tau_j (D^2u)\hat\tau_j|\Bigr\}(x_j,t_j)\ge c\, \delta^{-\beta p}(x_j)+j.
\end{equation} 
In view of \eqref{estBernstein1b}-\eqref{estBernstein2b} and of parabolic estimates, we have
$\delta(x_j)\to 0$.
Set
$$z_j:=P(x_j),$$
so that
$$\lambda_j:=\delta(x_j)=|x_j-z_j|.$$
After extracting a subsequence, we have $z_j\to a\in \partial\Omega$,
and we may assume without loss of generality that $a=0$ and $\nu_0=e_n$. 
Hence $x_j\to 0$, and
\be\label{defxij}
\xi_j:={x_j-z_j\over |x_j-z_j|}\to e_n.
\ee
We may also assume that
\be\label{HypContradtau}
\tau_j\to \tau_\infty,\quad \hat\tau_j\to \hat\tau_\infty,\quad\hbox{
with $\tau_\infty,\hat\tau_\infty \perp e_n$ and $|\tau_\infty|=1$, $|\hat\tau_\infty|=1$.}
\ee
Next we rescale $v$ by setting:
$$v_j(y,s):=\lambda_j^{\beta-1} u(z_j+\lambda_j y,t_j+\lambda_j^2 s),
\quad\hbox{ in $\Omega_j\times (-t_j\lambda_j^{-2},0]$},$$
where $\Omega_j:=\lambda_j^{-1}(\Omega-z_j)$
and $\lambda_j=\delta(x_j)\to 0$.
We have
$$\begin{aligned}
\nabla v_j(y,s)&=\lambda_j^\beta \nabla u(z_j+\lambda_j y,t_j+\lambda_j^2 s),\\
\Delta v_j(y,s)&=\lambda_j^{\beta+1} \Delta u(z_j+\lambda_j y,t_j+\lambda_j^2 s),\\
\partial_s v_j(y,s)&=\lambda_j^{\beta+1} u_t(z_j+\lambda_j y,t_j+\lambda_j^2 s).
\end{aligned}$$
Since $\beta+1=\beta p$, it follows that $v_j$ satisfies
$$\partial_s v_j-\Delta v_j=|\nabla v_j|^p
\quad\hbox{ in $\Omega_j\times (-t_j\lambda_j^{-2},0]$}.$$
Also $\Omega_j$ converges to the half-space $\{y_n>0\}$
and $t_j\lambda_j^{-2}\ge T/(2\lambda_j^{2})\to\infty$, as $j\to \infty$.

We now proceed to prove a suitable local compactness property of the sequence $(v_j)_j$
by making use of the Bernstein estimates \eqref{estBernstein1b}-\eqref{estBernstein2b}.
To this end we need to convert these estimates in terms of the variable $y$.
Since $\nu_0=e_n$, we first note that
\be\label{controldistj0}
\sup\Bigl\{\frac{|x_n-x'_n|}{|x-x'|},\ x, x'\in\partial\Omega, \, x\ne x',\, |x|+|x'|\le\eta\Bigr\} \to 0\quad\hbox{ as $\eta\to 0$}.
\ee
Since $\partial\Omega_j=\lambda_j^{-1}(\partial\Omega-z_j)$ and $z_j\in\partial\Omega$, $z_j\to 0$, $\lambda_j\to 0$,
it follows from \eqref{controldistj0}, applied with $x=z_j+\lambda_j\xi$ and $x'=z_j$, that, for any $R>0$,
\be\label{controldistj1}
\sigma_j(R):=\sup\bigl\{|y_n|;
\ y\in\partial\Omega_j, \, |y|\le 2R\bigr\} \to 0\quad\hbox{ as $j\to\infty$}.
\ee
Also, there exists $j_0(R)$ such that for any $j\ge j_0(R)$ and any $y\in \Omega_j\cap B_R$, the projection of $y$ onto $\partial\Omega_j$,
denoted by $\bar y^j$, is well defined and satisfies
\be\label{controldistj2}
y_n-\bar y^j_n\le |y-\bar y^j|={\rm dist}(y,\partial\Omega_j)\le 2(y_n-\bar y^j_n).
\ee
Consequently, for any $y\in \Omega_j\cap B_R$, since ${\rm dist}(y,\partial\Omega_j)\le |y|\le R$, we have
$$|\bar y^j|\le |\bar y^j-y|+|y|\le 2R.$$
 Therefore, \eqref{controldistj1} and \eqref{controldistj2} give
\be\label{controldistj}
y_n-\sigma_j(R)\le {\rm dist}(y,\partial\Omega_j)\le 2y_n+2\sigma_j(R),\quad y\in \Omega_j\cap B_R.
\ee
Fix any $R, \eta$, with $0<\eta<1<R$.
By \eqref{controldistj0} and \eqref{controldistj}, there exists $j_1(R,\eta)$ such that
\be\label{controldistyn0}
D_{R,\eta}:=\bigl\{y\in B_R;\,y_n>\eta\bigr\}\subset \Omega_j
\quad\hbox{ for all $j\ge j_1(R,\eta)$}
\ee
and
\be\label{controldistyn}
\frac12 y_n\le {\rm dist}(y,\partial\Omega_j)\le 3y_n
\quad\hbox{ for all $y\in D_{R,\eta}$ and all $j\ge j_1(R,\eta)$.}
\ee
Since
\be\label{controldistyn2}
{\rm dist}(z_j+\rho_j y,\partial\Omega)=\rho_j \,{\rm dist}(y,\partial\Omega_j),
\ee
we deduce from
 \eqref{estBernstein1b}, \eqref{estBernstein2b} and \eqref{controldistyn}
that, for all $j\ge j_1(R,\eta)$ 
and all $(y,s)\in D_{R,\eta}\times (-t_j\lambda_j^{-2},0]$, 
\be\label{BernsteinEst1}
|\nabla v_j(y,s)|
=\rho_j^\beta|\nabla u_j(z_j+\rho_j y, t_j+\lambda_j^2 s)|
\le C\rho_j^\beta \bigl[{\rm dist}(z_j+\rho_j y,\partial\Omega)\bigr]^{-\beta} \\
\le 2^\beta Cy_n^{-\beta}
\ee
and
\be\label{BernsteinEst2}
\begin{aligned}
|v_j(y,s)| &=\rho_j^{1-\beta}|u_j(z_j+\rho_j y,t_j+\lambda_j^2 s)|\\
&\le C\rho_j^{1-\beta} \bigl[{\rm dist}(z_j+\rho_j y,\partial\Omega)\bigr]^{1-\beta}
\le 3^{1-\beta}Cy_n^{1-\beta}.
\end{aligned}
\ee
By interior parabolic estimates, it follows that 
 the sequence $(v_j)_j$ is precompact in $C^{2,1}(Q_{R,\eta})$,
where $Q_{R,\eta}=D_{R,\eta}\times [-R,0]$. 
By a diagonal procedure, we deduce that some subsequence of $(v_j)_j$, not
relabeled, converges
in each $C^{2,1}(Q_{R,\eta})$ to a classical solution $V(y,s)\in C^{2,1}(\Rnp\times (-\infty,0])$ of
$$V_s-\Delta V=|\nabla V|^p,\quad \mbox{in }\Rnp\times (-\infty,0].$$
On the other hand, \eqref{boundutsol} yields
$$|\partial_s v_j(y,s)|\le M\lambda_j^{\beta+1},$$
so that we actually have $V_s=0$, hence $V=V(y)$.
Moreover, \eqref{BernsteinEst2} guarantees that
$$|V(y)|\le  3^{1-\beta}Cy_n^{1-\beta},\quad y_n\ge0.$$
 Consequently, $V$ extends to a function $V\in C(\overline{\Rnp})$, with $V=0$ on $\partial\Rnp$.
It follows from Theorem~\ref{LiouvilleThm} that either $V=0$ or $V(y)=V_\alpha(y_n)$ for some $\alpha\ge 0$.

Now, in case \eqref{HypContradDelta0a} holds, we have
$$\bigl|\nabla v_j(\xi_j,0)\bigr|\ge (1+c)d_p.$$
We then deduce from \eqref{defxij} that
$$\bigl|\nabla V(e_n)\bigr|\ge (1+c)d_p.$$
But this contradicts the fact that either $V=0$ or 
$\bigl|\nabla V(e_n)\bigr|=V'_\alpha(1)=d_p(\alpha+1)^{-\beta}\le d_p$.
We thus conclude that \eqref{est1thmBernstein2} is true.
Estimate \eqref{est2thmBernstein2} then follows by integration in the normal direction and this completes the proof of 
Theorem~\ref{Bernstein2}.

Finally, in case \eqref{HypContradDelta0} holds, we have
$$\bigl|\nabla v_j(\xi_j,0)\cdot\tau_j\bigr|^p+
\bigl|\tau_j \bigl(D^2v_j(\xi_j,0)\bigr)\nu(x_j)\bigr|+
\bigl|\tau_j \bigl(D^2v_j(\xi_j,0)\bigr)\hat\tau_j \bigr|\ge c.$$
 Hence, using \eqref{defxij}, \eqref{HypContradtau}, we get
$$\bigl|\nabla V(e_n)\cdot\tau_\infty\bigr|^p+
\bigl|\tau_\infty (D^2V(e_n))e_n\bigr|+
\bigl|\tau_\infty (D^2V(e_n))\hat\tau_\infty)\bigr|\ge c.$$
This contradicts the fact that $V$ depends only on $y_n$
and completes the proof of Proposition~\ref{lem1}. \qed

\medskip

\begin{rem}\label{remext} {\rm
Let $\bar u$ be the global viscosity solution of \eqref{main_pb} (cf.~the paragraph before 
Theorem~\ref{thmLBC}) and assume that,  for some $\tau\in (T,\infty]$,
\be\label{HypNoloss}
u=0 \quad\hbox{ on $\partial\Omega\times[0,\tau]$}
\ee
i.e.,~$\bar u$ does not lose boundary conditions for $t\le \tau$.
Then, for any $\eps>0$, estimates \eqref{est1thmBernstein2} and \eqref{est2thmBernstein2}
in Theorem~\ref{Bernstein2} remain valid for $\bar u$ in $\Omega\times[0,\tau]$.

Indeed, we know that $\bar u$ satisfies $\sup_{\Omega\times [T/2,\infty)} |\bar u_t| <\infty$
and that the gradient estimate~\eqref{estBernstein1b} remains valid 
for $\bar u$ in $\Omega\times[0,\infty)$.
This follows respectively from \cite[Lemma~10.1]{PS_jmpa18} and from \cite[Theorem~3.1]{PS_jmpa18}, 
applied with $F(\xi)=|\xi|^p$ and $\theta=(p-1)/p$.
Now, by integrating~\eqref{estBernstein1b} in the normal direction and using \eqref{HypNoloss},
we see that estimate~\eqref{estBernstein2b} remains valid 
for $\bar u$ in $\Omega\times[0,\tau)$.
The proof of Theorem~\ref{Bernstein2} then applies without changes.}
\end{rem}

\section{Proof of Theorem~\ref{thmEllipt}}
The proof is similar to that of Theorem~\ref{Bernstein2},
based on a rescaling argument, combined with the elliptic Bernstein estimates \cite{Lions85}
and the Liouville-type Theorem~\ref{LiouvilleThm}.
\smallskip

Assume for contradiction that there exist $d>d_p$ and sequences $\{f_j\}, \{u_j\}, \{x_j\}$,
with $\|f_j\|_{W^{1,\infty}}\le M$, such that
\be\label{hypcontradEll}
\delta(x_j)\to 0 \quad\hbox{and}\quad \delta^\beta(x_j)|\nabla u_j(x_j)|\ge d.
\ee
Let
$$z_j:=P(x_j),\quad \lambda_j:=\delta(x_j)=|x_j-z_j|.$$
By extracting a subsequence, we may assume without loss of generality that $z_j\to 0\in \partial\Omega$,
hence $x_j\to 0$, and
\be\label{convnu}
\nu_j:=\nu(z_j)\to e_n.
\ee
Set
$$v_j(y)=\lambda_j^{\beta-1}u_j(z_j+\lambda_j y).$$
We have
$$\nabla v_j(y)=\lambda_j^\beta \nabla u_j(z_j+\lambda_j y),\quad \Delta v_j(y)=\lambda_j^{\beta+1} \Delta u_j(z_j+\lambda_j y).$$
Since $\beta+1=\beta p$, it follows that $v_j$ satisfies
$$\Delta v_j+|\nabla v_j|^p=\lambda_j^{\beta+1}\bigl[\Delta u_j+|\nabla u_j|^p\bigr](z_j+\lambda_j y),$$
in $\Omega_j:=\lambda_j^{-1}(\Omega-z_j)$, hence
$$-\Delta v_j=|\nabla v_j|^p+\tilde f_j(y), \quad\hbox{where } \tilde f_j(y):=\lambda_j^{\beta+1}f_j(z_j+\lambda_j y)
\quad\hbox{ in $\Omega_j$.}$$
Also $\Omega_j$ converges to the half-space $\{y_n>0\}$ as $j\to \infty$.

On the other hand, it follows from the elliptic Bernstein estimate
\eqref{BernsteinPLL} that
\begin{align*}
|\nabla u_{j}(x)|&\le C\delta^{-\beta}(x)\quad\hbox{ in $\Omega$,}\\
|u_{j}(x)|&\le C\delta^{1-\beta}(x)\quad\hbox{ in $\Omega$,} 
\end{align*}
for some constant $C=C(n,p,M)>0$ independent of $j$.
Setting $Q_{R,\eps}=B_R\cap\{y_n>\eps\}$, arguing as in the proof of Proposition \ref{lem1} and using interior elliptic $L^q$ estimates,
we can find a subsequence of $v_j$, not relabeled, which converges
in $W^{2,q}(Q_{R,\eps})$ for each $q,R,\eps>0$ to a strong, hence classical, solution $V(y)\in C^2(\Rnp)\cap C(\overline{\Rnp})$ of
$$\begin{cases}
-\Delta u =|\nabla u|^p, &y\in\Rnp,
\\ u(y)=0, &y\in \partial\Rnp.
\end{cases}
$$
It follows from Theorem~\ref{LiouvilleThm} that
 either $V=0$ or $V(y)=U_\alpha(y_n)$ for some $\alpha\ge 0$.
In particular, $|\nabla V(y)|\le d_py_n^{-\beta}$ by \eqref{RefU0prime}. Since, by \eqref{convnu},
$$\lim_{j\to\infty}\delta^\beta(x_j) |\nabla u_j(x_j)|=\lim_{j\to\infty} |\nabla v_j(\nu_j)|=|\nabla V(e_n)|\le d_p.$$
 This contradicts \eqref{hypcontradEll}, being $d>d_p$.
\qed

\section{Proof of Theorem~\ref{ODEbehavior}}

The proof of Theorem~\ref{ODEbehavior} relies on Proposition \ref{lem1}, along with two other preliminary results.
The first one, Proposition~\ref{lemCentral}, which is rather technical,
is a nondegeneracy property for GBU points. It states that if the 
 singularity of $u_\nu$ is sufficiently weak 
 at small (but positive) distance from a given boundary point and in some time interval,
 then $|\nabla u|$ satisfies a uniform space-time bound near that point.

Proposition~\ref{lemCentral} is an improvement of \cite[Lemma 2.2]{LS10}, where the weak singularity assumption
had to be made in a whole space-time neighborhood of the boundary point (and not only at positive distance).
This improvement is crucial to our arguments and is made possible by
relying on estimate \eqref{ConclLemGrad} from Proposition~\ref{lem1}.
The proof amounts to showing that the weak singularity will propagate from
small finite distance up to the boundary.
Moreover, in the course of the proof of Theorem~\ref{ODEbehavior},
Proposition~\ref{lemCentral} will be applied to get uniform gradient estimates for
a sequence of suitably renormalized versions of $u$ in rescaled domains.
 Thus we need a local statement in more general domains under precise regularity
 assumptions on the boundary.

\begin{prop}\label{lemCentral}
Let $p>2$, let $\omega$ be a bounded domain of class $C^2$, $a\in\partial\omega$, $R>0$,
and set $D=B_{2R}(a)\cap \omega$. Assume that
\be\label{innersphere}
\hbox{$\omega$ satisfies an inner sphere condition of radius $R$ at each $b\in B_R(a)\cap \partial\omega$},
\ee
\be\label{distanceC2}
\hbox{the distance function $\delta(x)={\rm dist}(x,\partial\omega)$ is of class $C^2$ in $\overline D$}.
\ee
 Set $L:=\|\Delta \delta\|_{L^\infty(D)}$. 
Let $T_0\in\R$, $\theta>0$ and let $v\in C^{2,1}(\overline D\times[T_0-\theta,T_0))$ be a solution of
$$ 
\begin{cases}
v_t-\Delta v =|\nabla v|^p, & \mbox{in }D\times (T_0-\theta,T_0),\\
v=0, &\mbox{on }\big(B_R(a)\cap\partial\omega\big)\times (T_0-\theta,T_0),
\end{cases}
$$ 
satisfying
\be\label{boundut1central}
|v_t| \le M \quad \hbox{말n $D\times (T_0-\theta,T_0)$}
\ee
and the property that for any $\eps\in (0,1)$ there exists a constant $C_\eps>0$ such that
 for any vector field $\tau\perp \nu$,
\be\label{ConclLemGrad1central}
\Bigl[|\nabla v\cdot\tau|^p+|\tau (D^2v)\tau|\Bigr](x,t)\le \eps \delta^{-\beta p}(x)+C_\eps
\quad  \hbox{말n $D\times (T_0-\theta,T_0)$}.
\ee
For each $k\in(0,d_p)$ there exists $r_0>0$, depending only on $p,k,M,L$ and on the constants $C_\eps$, such that, if
 for some $r\in (0,\min\{r_0,R\}]$ and some $\sigma\in (0,r]$
\be\label{HypTypeILemSmall1central}
|v_\nu(b+r\nu_b,t)| \le kr^{-\beta},
\quad\hbox{ for all $b\in B_\sigma(a)\cap\partial\omega$ and all $t\in (T_0-\theta,T_0)$,}
\ee
 then
\be \label{ConclLemCentral}
|\nabla v|\le C(n,p,k,M,\sigma,\theta)
 \quad\hbox{ in $(B_{\sigma/8}(a)\cap\omega)\times (T_0-\theta/4,T_0)$}.
 \ee
\end{prop}

The proof of Proposition~\ref{lemCentral}, which is rather long and technical, is postponed to the appendix.
\medskip

Our last preliminary result, Proposition~\ref{lemCentral2}, gives an (optimal) lower bound on the final space profile
of $|u_\nu|$ in the normal direction
to a GBU point (however, the absolute value will be eventually removed in Theorem~\ref{thmprofile1}). It will be proved as a direct consequence of Proposition~\ref{lemCentral}.

\begin{prop}[] \label{lemCentral2}
Let $p>2$ and let $u_0\in X$ be such that $T:=T(u_0)<\infty$.
If $a$ is a GBU point of $u$
(i.e., $\limsup_{t\to T^-,\, x\to a} |\nabla u(x,t)|=\infty$),
 then
\begin{equation}
\label{HypTypeILemSmall2}
\liminf_{r\to 0}r^\beta |u_\nu(a+r\nu_a,T)| \ge d_p. 
\end{equation}
\end{prop}

\medskip

{\it Proof of Proposition \ref{lemCentral2}.}
Assume for contradiction that \eqref{HypTypeILemSmall2} fails.  Hence,
\be\label{liminfunusmall}
\liminf_{r\to 0}r^\beta |u_\nu(a+r\nu_a,T)|<k 
\ee
for some $k\in (0,d_p)$.
Take $R\in (0,\delta_0/2)$ such that $\Omega$ satisfies an inner sphere condition of radius $R$ at each point of $\partial\Omega$.
In view of \eqref{boundutsol} and Proposition~\ref{lem1}, we may apply Proposition~\ref{lemCentral} to $v=u$
with $\omega=\Omega$, $T_0=T$, $\theta=T/2$.
Let $r_0$ be given by Proposition~\ref{lemCentral} for the above value of $k$.
By \eqref{liminfunusmall},  there exists $r\in (0,\min\{r_0,\,R\}]$ such that
$$|u_\nu(a+r\nu_a, T)|<kr^{-\beta}.$$  
Since $u\in C^{2,1}(\Omega\times (0,T])$, by continuity, there exist $\sigma\in (0,r)$ and $\theta\in(0,T/2)$ so small that
$$|u_\nu(b+r\nu_b,t)|<kr^{-\beta}\quad\hbox{ for all $b\in B_\sigma(a)\cap\partial\Omega$ and all $t\in [T-\theta,T)$.}$$
It follows from Proposition~\ref{lemCentral}, that $a$ is not a GBU point: a contradiction.
\qed

\medskip

We are now in a position to give the proof of Theorem~\ref{ODEbehavior},
by combining Propositions~\ref{lem1}-\ref{lemCentral2} and an appropriate rescaling argument.
As mentioned in Remark~\ref{rem1}, we shall adapt the strategy in \cite{MZ98} to our problem,
also using some simplifications from \cite{GuoSou2}.
\medskip

{\it Proof of Theorem~\ref{ODEbehavior}}.
To establish \eqref{boundnormal}, it suffices to prove \eqref{boundtangential}, 
since \eqref{boundnormal} then follows in view \eqref{main_pb} and \eqref{boundutsol}.
We shall actually prove estimate \eqref{boundtangential2} at the same time,
 as noted in Remark~\ref{rem1b}(c).

Assume that either \eqref{boundtangential} or \eqref{boundtangential2} fails.
Thus there exist $c_1>0$, a sequence of couples $(x_j,t_j)\in \Omega_{\delta_0}\times [T/2,T)$
and unit vectors $\tau_j,\hat\tau_j$ with $\tau_j$ and $\hat\tau_j \perp \nu(x_j)$, such that
\begin{equation}
\label{HypContradDelta2}
K_j:=\Bigl\{|\nabla u\cdot\tau_j|^p+|\tau_j (D^2u)\nu|+|\tau_j (D^2u)\hat\tau_j|\Bigr\}(x_j,t_j)
\ge c_1|u_\nu(x_j,t_j)|^p+j.
\end{equation}
Set
$$z_j:=P(x_j),$$
so that $\delta(x_j)=|x_j-z_j|$.
The proof will be done in several steps. 
\smallskip

{\bf Step 1.} {\it Nondegeneracy at points $z_j$.}
First, it follows from Proposition~\ref{lem1} and~\eqref{HypContradDelta2} that, for all $\eps>0$ there exists $C_\eps>0$ such that
$$c_1|u_\nu(x_j,t_j)|^p+j\le K_j\le \eps \delta^{-\beta p}(x_j)+C_\eps.$$
Therefore, $t_j\to T$, $\delta(x_j)\to 0$ and
\be\label{HypSmallXj}
|u_\nu(x_j,t_j)|\le \eps_j \delta^{-\beta}(x_j),
\ee
with $\eps_j\to 0$ as $j\to\infty$.
After extracting a subsequence, we have $z_j\to a\in \partial\Omega$,
and we may assume without loss of generality that $a=0$, hence $x_j\to 0$, and that 
\be\label{HypCvNu}
\nu_j:=\nu(z_j)\to e_n.
\ee
Note that $0$ is in particular a GBU point (i.e., $\limsup_{t\to T^-,\, x\to 0} |\nabla u(t,x)|=\infty$),
since otherwise by parabolic regularity, $K_j$ would be bounded.

We claim that there exists a subsequence of $\{(x_j,t_j)\}_j$, not relabeled, and a sequence $\rho_j\to 0$ such that
\begin{equation}
\label{ClaimKappa}
\rho_j>\delta(x_j)\quad\hbox{and}\quad\rho_j^{\beta} |u_\nu(z_j+\rho_j\nu_j, t_j)|=\frac{d_p}{2}. 
\end{equation}
To prove the claim,
in view of \eqref{HypSmallXj}, by continuity, it suffices to show
that, for each $\rho>0$ and each $j_0\ge 1$ there exist $j\ge j_0$ and
$s\in(\delta(x_j),\rho)$ such that
$$ s^{\beta} |u_\nu(z_j+s\nu_j, t_j)|\ge \frac{d_p}{2}.$$ 
If this were false, then there would exist $\rho>0$ and $j_0\ge 1$ such that,
for all $j\ge j_0$ and $s\in (\delta(x_j),\rho)$, $s^{\beta} u_\nu(z_j+s\nu_j, t_j)< d_p/2$. 
Hence, letting $j\to\infty$,
$$s^{\beta} |u_\nu(se_n,T)|\le \frac{d_p}{2} 
\quad\hbox{ for all $s\in (0,\rho)$}.$$
Applying Proposition~\ref{lemCentral2}, we would deduce that $0$ 
is not a GBU point, which is a contradiction. This proves the claim.

\smallskip

{\bf Step 2.} {\it Rescaling and convergence to a one dimensional profile.}
We rescale similarly as in the proof of Proposition \ref{lem1},
but now taking $\rho_j$ as rescaling parameter.
Namely, we set:
\be\label{eqnvj}
v_j(y,s):=\rho_j^{\beta-1} u(z_j+\rho_j y,t_j+\rho_j^2 s)
\quad\hbox{ in $\Omega_j\times (-t_j\rho_j^{-2},0]$},
\ee
where $\Omega_j:=\rho_j^{-1}(\Omega-z_j)$.
The function $v_j$ satisfies
\be\label{main_pb_vj}\begin{cases}
\partial_s v_j-\Delta v_j=|\nabla v_j|^p&\quad\hbox{ in $\Omega_j\times (-t_j\rho_j^{-2},0]$}\\
v_j=0&\quad\hbox{on $\partial\Omega_j\times (-t_j\rho_j^{-2},0]$}\\
\end{cases}
\ee
and $\Omega_j$ converges to the half-space $\{y_n>0\}$ as $j\to \infty$.
Let $D_{R,\eta}:=\bigl\{y\in B_R;\,y_n>\eta\bigr\}$ and $Q_{R,\eta}= D_{R,\eta}\times [-R,0]$ for $0<\eta<1<R$.
Arguing exactly as in the proof of Proposition~\ref{lem1},
we can find a subsequence of $(v_j)_j$, not relabeled, which converges
in each $C^{2,1}(Q_{R,\eta})$ to a classical solution $w(y,s)\in C^{2,1}(\Rnp\times (-\infty,0])$ of
$$w_s-\Delta w=|\nabla w|^p
\quad\hbox{ in $\Rnp\times (-\infty,0]$.}$$
Moreover, by \eqref{boundutsol} we have
$$|\partial_s v_j(y,s)|\le M\rho_j^{\beta+1}.$$
 Hence $w_s\equiv 0$.
On the other hand, by (the analogues of) \eqref{controldistyn0}-\eqref{controldistyn2}, we have
$$\delta(z_j+\rho_j y) \ge \frac12 \rho_j y_n
\quad\hbox{ for all $j\ge j_1(R,\eta)$ and all $y\in D_{R,\eta}\subset\Omega_j$.}$$
 By Proposition~\ref{lem1}, for any $\eps>0$
 there exists $C_\eps>0$ such
that, for any unit vector field $\tau\perp\nu$,
 \begin{align*}
\bigl|\nabla v_j(y,s)\cdot\tau(z_j+\rho_j y)\bigr|
&=\rho_j^\beta \bigl|[\nabla  u\cdot\tau](z_j+\rho_j y,t_j+\rho_j^2 s)\bigr| \\
&\le \rho_j^\beta\bigl[\eps \delta^{-\beta}(z_j+\rho_j y)+C_\eps\bigr]\le 2^{\beta}\eps y_n^{-\beta}+C_\eps\rho_j^\beta.
\end{align*}
In view of \eqref{HypCvNu}, we deduce that $\partial_{y_1} w=\dots=\partial_{y_{n-1}} w\equiv 0$.
Therefore $w=w(y_n)$ and
the function $w$ solves
\be\label{solODE0}
-w''=|w'|^p,\quad y_n>0.  
\ee
Next, using $\nu_j\equiv\nu(z_j)=\nu(z_j+\rho_j \nu_j)$, property \eqref{ClaimKappa} yields
$$\bigl|\nabla v_j(\nu_j,0)\cdot\nu(z_j)\bigr|
=\rho_j^\beta \bigl|[\nabla u\cdot\nu](z_j+\rho_j \nu_j,t_j)\bigr|
=\rho_j^{\beta} |\partial_\nu u(z_j+\rho_j\nu_j,t_j)|=\frac{d_p}{2},
$$
so that $|w'(1)|=d_p/2$. But then necessarily $w'(1)=d_p/2$, 
since any solution of \eqref{solODE0}, with $w'(1)<0$, ceases to exist after some finite $y_n>1$.
Integrating \eqref{solODE0}, we thus obtain
$$ 
w'(y_n)=\bigl[\bigl(d_p/2\bigr)^{1-p}+(p-1)(y_n-1)\bigr]^{-\beta},\quad y_n>0. 
$$ 
On the other hand, letting
$$r_j:=\frac{\delta(x_j)}{\rho_j}\in (0,1),$$
we deduce from \eqref{HypContradDelta2} that 
\be
\label{HypContradDelta2b}
\begin{aligned}
&\bigl|\nabla v_j(r_j\nu_j,0)\cdot\tau_j\bigr|^p+\bigl|\tau_j\bigl(D^2v_j(r_j\nu_j,0)\bigr)\nu(x_j)\bigr|
+\bigl|\tau_j\bigl(D^2v_j(r_j\nu_j,0)\bigr)\hat\tau_j\bigr|\\
&\ \ =\rho_j^{\beta p}\bigl\{\bigl|\nabla u(x_j,t_j)\cdot \tau_j\bigr|^p+\bigl| \tau_j\bigl(D^2u(x_j,t_j)\bigr)\nu(x_j)\bigr|+
\bigl|{\hat\tau_j}\bigl(D^2u(x_j,t_j)\bigr)\hat\tau_j\bigr|\bigr\}\\
&\ \ \ge c_1 \rho_j^{\beta p}|\nabla u(x_j,t_j)\cdot\nu(x_j)|^p
=c_1|\nabla v_j(r_j\nu_j,0)\cdot\nu(x_j)|^p.
\end{aligned}
\ee
Since $w$ depends only on $y_n$, we expect to reach a contradiction with \eqref{HypContradDelta2b}.
However, since $r_j$ may approach $0$
and the convergence obtained so far is only valid locally for $y_n>0$ bounded away from $0$,
 we need to extend the convergence near $y_n=0$.
To this end, in the next step, we shall apply Proposition~\ref{lemCentral} to get a priori estimates of $|\nabla v_j|$ near the boundary.

\smallskip

{\bf Step 3.} {\it Uniform regularity of rescaled solutions and conclusion.}
 Put
$$A_\eta:=\{y\in\R^n;\, |y|\le 4,\, y_n\ge \eta\},\quad 0<\eta<1.$$
First, since
$$|\nabla w|=w'\le K_0:=\bigl[\bigl(d_p/2\bigr)^{1-p}-(p-1)\bigr]^{-\beta}
=\bigl(2^{p-1}-1)^{-\beta}d_p,$$
we deduce from Step~2 that
for any $\eta\in (0,1)$, there exists $j_1(\eta)\ge1$ such that
\be\label{boundvj}
|\nabla v_j|\le K_0+1
\quad\hbox{ in $A_\eta\times [-4,0],\ \ j\ge j_1(\eta)$.}
\ee

We want to show that the assumptions of Proposition~\ref{lemCentral} with  $\omega=\Omega_j=\rho_j^{-1}(\Omega-z_j)$
are satisfied, 
uniformly with respect to $j$ large.
It is easy to see that \eqref{innersphere} is satisfied for some $R\in (0,1)$ independent of $j$,
and moreover \eqref{distanceC2} is true with $D=D_j:=B_{2R}\cap \Omega_j$
(see~e.g.~\cite{CrM} for regularity properties of the function distance to the boundary).
Set
$$\delta_j(y)={\rm dist}(y,\partial\Omega_j), \quad \delta(x)={\rm dist}(x,\partial\Omega),$$
and observe that the normal vector field to $\partial\Omega_j$ is given by
$$\nu^j(y):=\nu(z_j+\rho_j y),$$
which is well defined in $D_j$.
Since $\delta_j(y)= \rho_j^{-1}\delta(z_j+\rho_j y)$, we have
$$|\Delta \delta_j(y)|=\rho_j|\Delta\delta(z_j+\rho_j y)|\le \rho_j\|\Delta \delta\|_{L^\infty(\Omega_{\delta_0})}
\le \|\Delta \delta\|_{L^\infty(\Omega_{\delta_0})},\quad y\in D_j.$$
Also, it follows from \eqref{boundutsol} that
$$|\partial_sv_j(y,s)|=\rho_j^{\beta+1} |u_t(z_j+\rho_j y,t_j+\rho_j^2 s)|\le M\rho_j^{\beta+1}\le M
\quad\hbox{ in $\Omega_j\times [-4,0],$}$$
where $M$ is independent of $j$.
Moreover, by Proposition~\ref{lem1} for any $\eps\in (0,1)$  there exists $C_\eps>0$, independent of $j$, such that
for any vector field $\tau(y)\perp \nu^j(y)=\nu(z_j+\rho_j y)$
$$\begin{aligned}
&\Bigl[|\nabla v_j\cdot\tau|^p+|\tau (D^2v_j)\tau|\Bigr](y,s)\\
&\qquad  =\rho_j^{\beta p}\Bigl\{\bigl|\nabla u(z_j+\rho_j y,t_j+\rho_j^2 s)\cdot\tau(y)\bigr|^p
+\bigl|\tau(y)\bigl(D^2u(z_j+\rho_j y,t_j+\rho_j^2 s)\bigr)\tau(y)\bigr|\Bigr\}\\
&\qquad\le\eps\rho_j^{\beta p}\,{\rm dist}^{-\beta p}(z_j+\rho_j y,\partial\Omega)+C_\eps \rho_j^{\beta p}\\
&\qquad=\eps\,{\rm dist}^{-\beta p}(y,\partial\Omega_j)+C_\eps \rho_j^{\beta p}
\le \eps\,{\rm dist}^{-\beta p}(y,\partial\Omega_j)+C_\eps
\quad\hbox{ in $D_j\times [-4,0]$}.
\end{aligned}$$
Now let $r_0=r_0(p,k,M,L,C_\eps)$ be given by Proposition~\ref{lemCentral} with $k=d_p/2$ ,
and choose $r\in (0,\min\{r_0,\,R\}]$ such that 
\be\label{boundvj2}
r^\beta(K_0+1)<d_p/2.
\ee
By the proof of \eqref{controldistj}, we have
$$
 {\rm dist}(y,\partial\Omega_j)\le 2y_n+2\sigma_j,\quad y\in \Omega_j\cap B_1,
$$
where
$\sigma_j:=\sup\bigl\{|\xi_n|;
\ \xi\in\partial\Omega_j, \, |\xi|\le 2\bigr\} \to 0$ as $j\to\infty$.
 Therefore,
$$\bigl\{b+r\nu^j_b;\, b\in\partial\Omega_j,\, |b|< r\bigr\}\subset A_{r/4}$$
 for all $j\ge j_2(r)$ large enough
which, along with \eqref{boundvj} and \eqref{boundvj2}, implies
$$
\Bigl|\partial_{\nu^j} v_j(b+r\nu_b,s)\Bigr| \le kr^{-\beta},
\quad\hbox{ for all $b\in B_r\cap\partial\Omega_j$, $s\in(-4,0]$, $j\ge j_2$.}
$$
We then deduce from Proposition~\ref{lemCentral} that
\be\label{controlnablaSigma}
|\nabla v_j|\le C(n,p,M,r)
 \quad\hbox{ in $(B_{r/8}\cap\Omega_j)\times (-2,0]$}.
 \ee
Writing $y=(y',y_n)$ and setting $\Sigma_j=\Omega_j\cap \{y;\ |y'|<r/16,\, |y|\le 3\}$,
we infer that, for all $j\ge j_3$ large enough,
\be\label{controlnablaSigma2}
|\nabla v_j|\le C'(n,p,M,r)
 \quad\hbox{ in $\Sigma_j\times (-2,0]$}
 \ee
(consider the cases $y_n\le r/16$ and $y_n> r/16$ and use \eqref{controlnablaSigma} and \eqref{boundvj}, respectively).

Next note that, for $j$ large enough, $\partial\Omega_j\cap B_4$ can be written as the graph
$\{y_n=\psi_j(y')\}$ of a $C^2$ function $\psi_j$,
with a uniform $C^2$-bound with respect to $j$.
Going back to \eqref{main_pb_vj}, 
we may thus apply parabolic interior-boundary estimates uniformly in $j$,
to deduce from \eqref{controlnablaSigma2} that
$$\sup_{j}\|v_j\|_{C^{2+\gamma,1+(\gamma/2)}(\Sigma_j'\times (-1,0])}<\infty,\quad\hbox{ where }
\Sigma_j'=\Omega_j\cap \{y;\ |y'|<r/32,\, |y|\le 2\},$$
for some $\gamma\in (0,1)$.
Passing to a subsequence and recalling that $\nu_j\to e_n$,
we may also assume that $\tau_j\to\tau_\infty$, $\tilde\tau_j\to \tilde\tau_\infty$,
with $|\tau_\infty|=|\tilde\tau_\infty|=1$, $\tau_\infty,\tilde\tau_\infty\perp e_n$, and that $r_j\to \ell\in [0,1]$.
In view of \eqref{HypContradDelta2b}, this implies
$$\begin{aligned}
\bigl|\nabla w(\ell e_n,0)\cdot\tau_\infty\bigr|^p
&+\bigl|\tau_\infty (D^2w(\ell e_n,0))e_n\bigr|
+\bigl|\tau_\infty (D^2w(\ell e_n,0))\hat\tau_\infty\bigr|\\
& \ge c_1|\partial_{y_n} w(\ell e_n,0)|^p=c_1|w'(\ell)|^p>0.
\end{aligned}$$
This contradicts the fact that $w=w(y_n)$ and the proof
of \eqref{boundtangential} and \eqref{boundtangential2} is completed. 
\smallskip

{\bf Step 4.} {\it Proof of \eqref{boundinf}.} Denote by $\lambda_1>0$ the first eigenvalue of $-\Delta$ in $H^1_0(\Omega)$,
and by $\varphi_1$ the corresponding positive eigenfunction, normalized by $\max_{\overline\Omega}\varphi_1=1$.
Since $u_0\in X$, there exists $M_1>0$ such that $u_0\ge -M_1\varphi_1$ in $\overline\Omega$.
Observing that $\underline u(x,t):= -M_1 e^{-\lambda_1 t}\varphi_1$ is a subsolution of~\eqref{main_pb},
if follows from the comparison principle that $u\ge \underline u$ in $\Omega\times (0,T)$.
 In particular,
$$-u_\nu\le M_2:=M_1\max_{\partial\Omega}\partial_\nu\varphi_1\quad\hbox{ on $\partial\Omega\times (0,T)$.}$$
But, by \eqref{boundnormal}, we have
$$|u_{\nu\nu}|\le 2|u_\nu|^p+C_1
\quad\hbox{ in $\Omega_{\delta_0}\times [T/2,T]$.}$$
For each $t\in [T/2,T)$ and $a\in \partial\Omega$, the function $\phi(s):=-u_\nu(a+s\nu_a,t)$ thus satisfies
$\phi'\le 2|\phi|^p+C_1$ in $[0,\delta_0]$, with $\phi(0)\le M_2$.
It then follows from standard ODE arguments that
\be\label{controlphi1}
\phi(s)\le M_2+1\quad\hbox{ for all $s\in [0,s_0]$,}
\ee
where $s_0\in (0,\delta_0)$ depends only on $p,M_2,C_1$.
But on the other hand, as a consequence of~\eqref{estBernstein1b}, we have
\be\label{controlphi2}
\phi(s)\le |\nabla u(a+s\nu_a)|\le Cs_0^{-\beta} \quad\hbox{ for all $s\in (s_0,\delta_0]$}.
\ee
Combining \eqref{controlphi1} and \eqref{controlphi2}, we obtain
$\displaystyle\inf_{\Omega_{\delta_0}\times (T/2,T)}u_\nu >-\infty$.
 This implies \eqref{boundinf}, owing to the $C^1$ spatial regularity of $u$ on $\overline\Omega\times [0,T/2]$.
 \qed

\section{Proof of Theorems~\ref{thmprofile1} and \ref{thmLBC}}

{\it Proof of Theorem~\ref{thmprofile1}.}
By Proposition~\ref{lemCentral2} and \eqref{boundinf} in Theorem~\ref{ODEbehavior}, we have
$$\liminf_{s\to 0}s^\beta \nabla u(a+s\nu_a,T)\cdot\nu_a \ge d_p.$$
Assertion (i) then follows from \eqref{est1thmBernstein2} in Theorem~\ref{Bernstein2},
being $\eps$ arbitrary.

\smallskip

Let us next prove assertion (iii).
For fixed $b\in \partial\Omega$ and $t\in [T/2,T)$, we set
$$\phi(r):=u_\nu(b+r\nu_b,t),\quad 0<r<\delta_0,$$
Since $\phi'(r)=u_{\nu\nu}(b+r\nu_b,t)$,
estimate \eqref{boundnormal} guarantees that
$$|\phi'|\le 2|\phi|^p+C_1\quad\mbox{in } (0,\,\delta_0),$$
for some constant $C_1=C_1(u_0)$.
For given $M>0$, take $\eps=\eps(M)\in (0,\delta_0)$ such that $[2(2M)^p+C_1]\eps<M$.
A standard argument shows that, if $|\phi(r)|\le M$
  at some $r\in (0,\eps)$, then $|\phi(s)|\le 2M$ for all
$s\in(0,r+\eps]$.

Now assume for contradiction that $u_\nu(x_j,t_j)<M$ for some $M>0$ and some sequences $x_j\to a$, $t_j\to T$. 
By \eqref{boundinf}, we thus have $|u_\nu(x_j,t_j)|<M$ for a possibly larger $M$. 
We may write $x_j=b_j+r_j\nu_{b_j}$ with $b_j=P(x_j)$ and $r_j=\delta(x_j)\to 0$.
Since $r_j<\eps(M)$ for all $j$ large, we deduce from the previous paragraph that
$$|u_\nu(b_j+s\nu_{b_j},t_j)|\le 2M \quad\hbox{ for all $s\in(0,r_j+\eps]$}.$$
Passing to the limit, we get
$$|u_\nu(a+s\nu_a,T)|\le 2M \quad\hbox{ for all $s\in(0,\eps]$}.$$ 
But this contradicts Proposition~\ref{lemCentral2} and so
assertion (iii) follows.
\smallskip

Let us now prove assertion (ii).
Fix $\eps\in (0,1)$. By the smoothness of $\partial\Omega$, we may find $\eta_\eps>0$ such that,
for all $x_0\in\partial\Omega\cap B_{\eta_\eps}(a)$,
 with $x_0\ne a$,
there exists a $C^1$ curve $\gamma: [0,1]\to \partial\Omega$
such that $\gamma(0)=a$ and $\gamma(1)=x_0$, and
there exists a unit vector $\tau\perp\nu_a$ such that
\be\label{regulgamma}
\sup_{s\in [0,1]} \left|\frac{\gamma'(s)}{|x_0-a|}-\tau \right|\le \eps.
\ee
Assertion (iii) 
gives
\be\label{inftylim}
\lim_{t\to T^-,\,x\to a} u_\nu(x,t)=\infty.
\ee
For $t\in (T/2,T)$ and $x_0\in\partial\Omega\cap B_{\eta_\eps}(a)$, we may then write 
$$u_\nu^{1-p}(x_0,t)-u_\nu^{1-p}(a,t)=-(p-1)\int_0^1 [u_\nu^{-p}\,\nabla u_\nu](\gamma(s),t)\cdot\gamma'(s)\,ds.$$
By \eqref{regulgamma}, we have
$$|\nabla u_\nu(\gamma(s),t)\cdot\gamma'(s)|
\le \Bigl(|u_{\nu\tau}(\gamma(s),t)|+\eps|\nabla u_\nu(\gamma(s),t)|\Bigr)|x_0-a|.$$
Since $|u_{\nu\tau}|\le \eps|u_\nu|^p+C_\eps $ and $|\nabla u_\nu|\le 2|u_{\nu}|^p+C$ in $\Omega_{\delta_0}\times [T/2,T]$ by Theorem~\ref{ODEbehavior}, it follows that
$$|\nabla u_\nu(\gamma(s),t)\cdot\gamma'(s)|
\le \Bigl(3\eps u_{\nu}^p(\gamma(s),t)+2C_\eps\Bigr)|x_0-a|.$$
 Hence, using \eqref{inftylim},
\be\label{inftylim2a}
u_\nu^{1-p}(x_0,t)-u_\nu^{1-p}(a,t)\le (p-1)|x_0-a|\int_0^1 [3\eps+2C_\eps u_\nu^{-p}](\gamma(s),t)\,ds \le 4(p-1)\eps|x_0-a|
\ee
for all $(x_0,t)\in [\partial\Omega\cap B_{\eta_\eps}(a)]\times (T-\eta_\eps,T)$, with $\eta_\eps>0$ possibly smaller. 
Now, if $x_0$ is not a GBU point, then $u_\nu(x_0,t)$ has a finite limit as $t\to T$. Letting $t\to T$ and using \eqref{inftylim} again, we obtain
\be\label{inftylim2}
u_\nu(x_0,T) \ge [4(p-1)\eps|x_0-a|]^{-\beta}
\ee
If $x_0$ is a GBU point, then $\lim_{t\to T^-}u_\nu(x_0,t)=\infty$, so that \eqref{inftylim2} remains true as well.
Since $\eps$ is arbirarily small, we have thus proved \eqref{moresingul}.

\smallskip

Finally, to check assertion (iv), we observe that 
\eqref{spacetime1} is a consequence of \eqref{boundnormal} 
and~\eqref{contnabla}.
As for \eqref{spacetime2}, it follows from \eqref{spacetime1} by integrating in the normal direction,
whereas \eqref{spacetime3} is a consequence of \eqref{inftylim2a}.
\qed

\medskip

{\it Proof of Theorem~\ref{thmLBC}.}
It suffices to prove assertion (ii), since assertion (i) will then clearly follow by exchanging the roles of $u_0$ and $v_0$.
Recall (see e.g.~\cite{PS_jmpa18}) that
$\bar u, \bar v\in C^{2,1}(\Omega\times (0,\infty))\cap C(\overline\Omega\times [0,\infty))$ and that,
since $v_0\ge u_0$, we have
$$T(v_0)\le T(u_0)\quad\hbox{ and \quad$\bar v\ge \bar u$ on $\overline\Omega\times (0,\infty)$. }$$
Pick $t_0\in (0,T(v_0))$. A consequence of Hopf's Lemma gives $v(\cdot,t_0)\ge
\lambda u(\cdot,t_0)$
for some $\lambda>1$.
Assume for contradiction that
\be\label{LBCcontrad}
\bar v=0 \quad\hbox{ on $\partial\Omega\times [0,T(u_0)]$.}
\ee
Then, by Theorem~\ref{Bernstein2} and Remark~\ref{remext}, we have
\be\label{LBClarge0}
\limsup_{\delta(x)\to 0} \delta^{\beta-1}(x) \, \bar v(x,T(u_0))\le c_p.
\ee
On the other hand, since $w:=\lambda u$ satisfies in $\Omega\times (0,T(u_0))$
$$w_t-\Delta w-|\nabla w|^p=\lambda(u_t-\Delta u-\lambda^{p-1}|\nabla u|^p)\le 0
 =\bar v_t-\Delta \bar v-|\nabla \bar v|^p$$
and $w=\bar v=0$ in $\partial\Omega\times (0,T(u_0))$, it follows from the comparison principle
(cf.~e.g.~\cite[Proposition 3.3]{PS_jmpa18}) that
$\bar v\ge \lambda u$ in $\Omega\times (t_0,T(u_0)]$.
Let $a\in \partial\Omega$ be a GBU point of $u$. By Theorem~\ref{thmprofile1}, we have
$$\lim_{s\to 0} s^{\beta-1} u(a+s\nu_a,T(u_0))=c_p.$$
Hence
$$\liminf_{s\to 0} s^{\beta-1} \bar v(a+s\nu_a,T(u_0))\ge \lambda c_p.$$
But this contradicts \eqref{LBClarge0}.
 Consequently, \eqref{LBCcontrad} cannot hold,
i.e.~$\bar v$ loses boundary conditions before $T(u_0)$.
 Hence in particular $T(v_0)<T(u_0)$. This proves assertion~(ii)
 and completes the proof.
\qed

\section{Appendix: proof of Proposition \ref{lemCentral}}

It is based on three lemmas (that we state together and will prove afterwards).
The first one, based on estimate \eqref{ConclLemGrad} from Proposition~\ref{lem1}, shows
that sufficiently weak singularity will propagate from
small finite distance up to the boundary.

\begin{lem}\label{lem2b} Let $p>2$, let $\omega\subset \R^n$ be an open set.
 Assume that the open line segment $(a,a+Re)\subset \omega$, where $a, e\in \R^n$, $|e|=1$, $R>0$.
Let $t_0<t_1$ and let $v\in C^{2,1}(\omega\times[t_0,t_1))$
be a solution of
\be\label{main_pbloc1}
v_t-\Delta v =|\nabla v|^p \quad \hbox{말n $\omega\times (t_0,t_1)$}
\ee
satisfying
\be\label{boundut1}
|v_t| \le M \quad \hbox{말n $\omega\times (t_0,t_1)$}
\ee
and with the property that for any $\eps>0$  there exists a constant $C_\eps>0$ such that for any unit vector $\tau\perp e$ 
\be\label{ConclLemGrad1}
\Bigl[|\nabla v\cdot\tau|^p+|\tau (D^2v)\tau|\Bigr](a+se,t)\le \eps  s^{-\beta p}+C_\eps,
\quad\hbox{ $s\in(0,R),\ t\in (t_0,t_1)$}.
\ee
For each $k\in(0,d_p)$, there exists $r_0>0$, depending only on $p,k,M$ and on the constants $C_\eps$, such that, if
\be\label{HypTypeILemSmall1a} 
[\nabla v\cdot e](a+re,t) \le kr^{-\beta},
\quad\hbox{ for some $r\in (0,\min\{r_0,R\}]$ and $t\in (t_0,t_1)$,}
\ee
 then
\be\label{ConclTypeILemSmall1}
[\nabla v\cdot e](a+se,t)\le k\,s^{-\beta}\quad\hbox{ for all $s\in (0,r]$.} 
\ee
\end{lem}

The second lemma provides regularization in time of the boundary derivative,
 provided that sufficiently weak singularity occurs. 
It is in the spirit of \cite[Lemma 2.2]{LS10}, and \cite[Lemma 4.1]{PS_imrn16} and, like these results, its proof relies on a barrier argument.
However, the statements there are not sufficient. In particular we need a uniform version
 for the purposes of the present paper.

\begin{lem}[] \label{lem0bApp}
Let $p>2$, let $\omega$ be a bounded domain of class $C^2$, $a\in\partial\omega$, $\rho>0$.
Set $D=B_\rho(a)\cap \omega$, assume that
the distance function $\delta(x)={\rm dist}(x,\partial\omega)$ is of class $C^2$ in $\overline D$, and set
\be\label{bounddelta}
L:=\|\Delta \delta\|_{L^\infty(D)}.
\ee
Let $t_0<t_1$ and let $v\in C^{2,1}(\omega\times[t_0,t_1))$ be a solution of
$$ 
\begin{cases}
v_t-\Delta v \le |\nabla v|^p, &  \mbox{in }\omega\times(t_0,t_1),\\ 
v=0, &  \mbox{on } \partial\omega\times(t_0,t_1).
\end{cases}
$$ 
For any $k\in (0,d_p)$, there exists $\rho_0=\rho_0(n,p,k,L)>0$ such that, if
$0<\rho\le \rho_0$ and
\be\label{HypTypeILemSmall1}
v \le (1-\beta)^{-1}k\delta^{1-\beta}
\quad\hbox{ in $D\times (t_0,t_1)$,}
\ee
then
$$
v_\nu \le c(n,p,k)\rho^{-\frac{1}{p-1}}\left(\frac{t_1-t_0+\rho^2}{t-t_0}\right)^{\frac{1}{p-2}}
\quad\hbox{ in $(B_{\rho/2}(a)\cap\partial\omega)\times (t_0,t_1)$}.
$$
\end{lem}

The last lemma provides a bound of the gradient near a boundary point
assuming a control of the normal derivative on the boundary.
The result and the proof (based on a local Bernstein-type argument)
are similar to \cite[Lemma 2.1]{LS10}, but again we require a more quantitative form.

\begin{lem}\label{lem:exclude blow-up points}
Let $p>2$, let $\omega$ be any domain of $\R^n$ and let $x_0\in\partial\omega$.
Let $R>0$, $t_0<t_1$ and $v\in C^{2,1}\bigl(
(\overline{B}_R(x_0)\cap\overline\omega)\times[t_0,t_1)\bigr)$ be a solution of
\be\label{main_pbloc}\begin{cases}
v_t-\Delta v =|\nabla v|^p, &
\mbox{in }(B_R(x_0)\cap\omega)\times(t_0,t_1),\\
v=0, &\mbox{on }(B_R(x_0)\cap\partial\omega)\times(t_0,t_1),
\end{cases}\ee
such that
\begin{equation}\label{cond:bounded gradient u}
    |\nabla v|\leq N\ \ \ {\rm on}\ \  (B_R(x_0)\cap\partial\omega)\times[t_0,t_1)
\end{equation}
and
\begin{equation}\label{cond:bounded ut}
    |v_t|\leq M\ \ \ {\rm in}\ \  (B_R(x_0)\cap\omega)\times[t_0,t_1),
\end{equation}
for some $M,N>0$.
Then there exists $C(n,p)>0$ such that
$$
    |\nabla v|\leq C(n,p)\Bigl[N+M^{1/p}+R^{-1/(p-1)}+(t-t_0)^{-1/2(p-1)}\Bigr]
$$
 in $\big(B_{R/2}(x_0)\cap\omega\big)\times(t_0,t_1)$.
\end{lem}

{\it Proof of Lemma \ref{lem2b}.}
Set
$$\phi(s):=[\nabla v\cdot e](a+se,t),\quad 0<s<R.$$
Then
\be\label{phiprime}
\phi'(s)=eD^2v(a+se,t)e.
\ee
Take any $\eta\in(0,1)$. By \eqref{ConclLemGrad1}, there exists $\delta_\eta>0$,
depending only on $\eta, p$ and on the constants~$C_\eps$, such that,
 for any unit vector $\tau\perp e$,
\be\label{ConclLemGrad1b}
|\nabla v(a+se,t)\cdot\tau|\le \eta  s^{-\beta}
\quad\hbox{and}\quad \bigl|\tau \bigl(D^2v\bigr)(a+se,t)\tau\bigr|\le \eta  s^{-\beta p},
\ee
for all $s\in (0,\min\{\delta_\eta,R\}]$ and $t\in (t_0,t_1)$.
Using the elementary inequality
\be\label{ineqab}
(X+Y)^{p/2}\le (1+\eta)X^{p/2}+c_0\eta^{1-(p/2)}Y^{p/2},\quad X,Y>0,
\ee
with $c_0=c_0(p)>0$, it follows that
$$|\nabla v(a+se,t)|^p\le \bigl|\phi^2(s)+\eta^2 s^{-2\beta}\bigr|^{p/2}
\le (1+\eta)|\phi|^p+c_0\eta^{1+\frac{p}{2}} s^{-\beta p}
$$
for all $s\in (0,\min\{\delta_\eta,R\}]$ and $t\in (t_0,t_1)$.
Using \eqref{main_pbloc1}, \eqref{boundut1}, \eqref{phiprime} and \eqref{ConclLemGrad1b}
and taking $\delta_\eta$ even smaller if necessary (depending only on $\eta,p,M$ and on the constants~$C_\eps$), we obtain
$$
\begin{aligned}\phi'(s)
&=e D^2v(a+re,t)e-\Delta v(a+se,t)+v_t(a+se,t)-|\nabla v(a+se,t)|^p \\
&\ge -(n-1)\eta s^{-\beta p}-M-(1+\eta)|\phi(s)|^p-c_0\eta^{1+\frac{p}{2}} s^{-\beta p} \\
&\ge -(1+\eta)|\phi(s)|^p-c_1\eta s^{-\beta p},\qquad 0<s<\delta_\eta,
\end{aligned}
$$
where $c_1=c_1(n,p)>0$.
Let now $\psi(s)=k s^{-\beta}$. Since $k<d_p$, hence $-\beta k+k^p<0$, 
we may choose $\eta>0$ sufficiently small
(depending only on $p,k,n$) so that
$-\beta k +(1+\eta)k^p+c_1\eta<0$. A simple computation gives
$$\psi'+(1+\eta)|\psi|^p+c_1\eta s^{-\beta p}=[-\beta k +(1+\eta)k^p+c_1\eta]s^{-\beta p}<0,\quad s>0.$$
Now choose $r_0=\delta_\eta$. Since $\phi(r)\le\psi(r)$ by \eqref{HypTypeILemSmall1a}, 
it follows from ODE comparison that $\phi(s)\le\psi(s)$ for all $s\in(0,r]$, i.e.~\eqref{ConclTypeILemSmall1}. \qed

\medskip

{\it Proof of Lemma \ref{lem0bApp}.}
 Fix a function $\Theta\in C^\infty_0(\mathbb R^n)$ with $\Theta(x)=1$ for $|x|\le 1/2$, $\Theta(x)=0$ for $|x|\ge 1$,
and $0\le\Theta\le 1$. Next set $\psi=\Theta^2$, $\tau=t_1-t_0$ and $\kappa=(1-\beta)^{-1}k$. We define
$$V=\kappa W,\qquad W(x,t)=\bigl(\delta(x)
+\varphi(x,t)\bigr)^{1-\beta}-\varphi^{1-\beta}(x,t),\quad \hbox{ for $(x,t)\in D\times (t_0,t_1),$}$$
with
$$\varphi(x,t)= \eta\rho \bigl[h(t)\psi(\rho^{-1}x)\bigr]^{1\over 1-\beta},\qquad h(t):=\tau^{-1}(t-t_0),$$
where $\eta\in(0,1)$ will be chosen suitably small below.
Note that $\varphi\in C^{2,1}(\R^n\times [t_0,t_1])$, $\varphi\ge 0$ and $0\le h\le 1$.
In this proof we denote by $C_1, C_2,\dots$ various positive constants depending only on $p$ and $n$ (through $\psi$).
We have
\begin{align}
\varphi&\le \eta\rho \label{derivphiA} \\
\varphi_t&\ge 0, \label{derivphiB}\\
\bigl(\varphi^{1-\beta}\bigr)_t
&=(\eta\rho)^{1-\beta}\tau^{-1}\psi(\rho^{-1}x)
\le (\eta\rho)^{1-\beta}\tau^{-1},   \label{derivphiC}\\ 
\bigl|\nabla(\varphi^{1-\beta})\bigr|
&=(\eta\rho)^{1-\beta}h(t)\rho^{-1}|\nabla\psi(\rho^{-1}x)|
\le C_1\eta^{1-\beta}\rho^{-\beta}, \\
\bigl|\Delta(\varphi^{1-\beta})\bigr|
&=(\eta\rho)^{1-\beta}h(t)\rho^{-2}|\Delta\psi(\rho^{-1}x)|
\le C_1\eta^{1-\beta}\rho^{-1-\beta}. \label{derivphiZ}
\end{align}
Morevover,
\be\label{nabladeltaphiA}
|\nabla\varphi|={\eta\rho\over 1-\beta}h^{1\over 1-\beta}(t)
\rho^{-1}\Bigl|\psi^{\beta\over 1-\beta}\nabla\psi\Bigr|(\rho^{-1}x)\le C_2\eta
\ee
and $\psi^{2\beta-1\over 1-\beta}|\nabla\psi|^2=4\Theta^{{2(2\beta-1)\over 1-\beta}+2}|\nabla\Theta|^2
=4\Theta^{2\beta\over 1-\beta}|\nabla\Theta|^2$ where $\psi\ne 0$. Hence,
\be\label{nabladeltaphiZ}
|\Delta\varphi|={\eta\rho\over 1-\beta}h^{1\over 1-\beta}(t)
\rho^{-2}\Bigl|\psi^{\beta\over 1-\beta}\Delta\psi+{\beta\over 1-\beta}\psi^{2\beta-1\over 1-\beta}|\nabla\psi|^2\Bigr|(\rho^{-1}x) \le C_2\eta\rho^{-1}.
\ee
Setting $Z=\bigl(\delta +\varphi\bigr)^{1-\beta}$, we compute
\begin{align*}
Z_t&=(1-\beta)(\delta +\varphi)^{-\beta}\varphi_t\ge 0, \\
\nabla Z&=(1-\beta)(\delta +\varphi)^{-\beta}\nabla(\delta +\varphi), \\
\Delta Z&=(1-\beta)(\delta +\varphi)^{-\beta}\Delta(\delta +\varphi)
-\beta(1-\beta)(\delta +\varphi)^{-\beta-1}|\nabla(\delta +\varphi)|^2.
\end{align*}
It follows that
$$W_t-\Delta W
\ge \beta(1-\beta)(\delta +\varphi)^{-\beta-1}|\nabla(\delta +\varphi)|^2
-(1-\beta)(\delta +\varphi)^{-\beta}\Delta(\delta +\varphi)
-\bigl(\varphi^{1-\beta}\bigr)_t+\Delta \bigl(\varphi^{1-\beta}\bigr).$$
Also, for $\eps\in (0,1)$ to be chosen below, by the elementary inequality
$$(X+Y)^p\le (1+\eps)X^p+C_3\eps^{1-p}Y^p,\quad X,Y>0,$$
we have
$$|\nabla W|^p\le (1+\eps) (1-\beta)^p(\delta +\varphi)^{-\beta p}\bigl|\nabla(\delta +\varphi)\bigr|^p
+ C_3\eps^{1-p}\bigl|\nabla(\varphi^{1-\beta})\bigr|^p.
$$
Therefore,
$$\begin{aligned}
&\kappa^{-1}\bigl(V_t-\Delta V-|\nabla V|^p\bigr) \\
&\quad\ge \beta(1-\beta)(\delta +\varphi)^{-\beta-1}|\nabla(\delta +\varphi)|^2
-(1-\beta)(\delta +\varphi)^{-\beta}\Delta(\delta +\varphi)
-\bigl(\varphi^{1-\beta}\bigr)_t+\Delta \bigl(\varphi^{1-\beta}\bigr)\\
&\qquad -\kappa^{p-1}\Bigl[(1+\eps) (1-\beta)^p(\delta +\varphi)^{-\beta p}\bigl|\nabla(\delta +\varphi)\bigr|^p
+C_3\eps^{1-p}\bigl|\nabla(\varphi^{1-\beta})\bigr|^p\Bigr].
\end{aligned}$$
Using $\beta p=\beta+1$, we obtain in the set $D\times (t_0,t_1)$
\begin{align*}
&\kappa^{-1}(\delta +\varphi)^{\beta p}\bigl(V_t-\Delta V-|\nabla V|^p\bigr) \\
&\quad\ge (1-\beta)|\nabla(\delta +\varphi)|^2
\Bigl[\beta-\kappa^{p-1}(1+\eps) (1-\beta)^{p-1}\bigl|\nabla(\delta +\varphi)\bigr|^{p-2}\Bigr] \\
&\qquad-(1-\beta)(\delta +\varphi)\Delta(\delta +\varphi)
+(\delta +\varphi)^{\beta p}\Bigl[-\bigl(\varphi^{1-\beta}\bigr)_t+
\Delta \bigl(\varphi^{1-\beta}\bigr)-C_3 \kappa^{p-1}\eps^{1-p}\bigl|\nabla(\varphi^{1-\beta})\bigr|^p\Bigr].
\end{align*}
On the other hand, owing to $|\nabla\delta|=1$ and \eqref{nabladeltaphiA}, we have
$$(1-C_2\eta)_+\le \bigl|\nabla(\delta +\varphi)\bigr|\le 1+C_2\eta.$$
Also, since $\beta>k^{p-1}$,
we may choose $\eps=\eps(p,n,k)\in(0,1)$ sufficiently small, so that
$$\beta-(1+\eps)(1+C_2\eps)^{p-2}k^{p-1}>0.$$
Assuming $\eta\in (0,\eps]$, using $(1-\beta)\kappa=k$, \eqref{derivphiA} \eqref{derivphiC}-\eqref{derivphiZ}, \eqref{nabladeltaphiZ}, $\delta\le 2\rho$ in $D$,
and recalling \eqref{bounddelta}, we derive
\begin{align*}
&\kappa^{-1}(\delta +\varphi)^{\beta p}\bigl(V_t-\Delta V-|\nabla V|^p\bigr) \\
&\quad\ge (1-\beta)(1-C_2\eta)_+^2
\Bigl[\beta-(1+\eps)(1+C_2\eta)^{p-2}k^{p-1}\Bigr]  -(\delta +\varphi)(L+C_2\eta\rho^{-1})\\
&\qquad-(\delta +\varphi)^{\beta p}\bigl[(\eta\rho)^{1-\beta}\tau^{-1}
+C_1\eta^{1-\beta}\rho^{-1-\beta}+C_3C_1^p
\kappa^{p-1}\eps^{1-p}\eta^{(1-\beta)p}\rho^{-\beta p}\bigr]\\ 
&\quad\ge (1-\beta)(1-C_2\eta)_+^2
\Bigl[\beta-(1+\eps)(1+C_2\eta)^{p-2}k^{p-1}\Bigr]  \\
&\qquad -(\eta+2)(L\rho+C_2\eta)-C_4\eta^{1-\beta}[(\eta+2)\rho]^{\beta p}\rho^{-\beta p}\bigl(\rho^2\tau^{-1}+1+\eps^{1-p}\eta^{p-2}\bigr)\\
&\quad\ge (1-\beta)(1-C_2\eta)_+^2
\Bigl[\beta-(1+\eps)(1+C_2\eta)^{p-2}k^{p-1}\Bigr] \\
&\qquad-3L\rho
-C_5\eta^{1-\beta}\bigl(\rho^2\tau^{-1}+1+\eps^{1-p}\bigr).
\end{align*}
 Now taking $\eps=\eps(p,n,k)\in(0,1)$ possibly smaller and $\rho_0=\rho_0(p,n,k,L)\in(0,1)$ sufficiently small,
we may assume that
 $$(1-\beta)(1-C_2\eps)_+^2\Bigl[\beta-(1+\eps)(1+C_2\eps)^{p-2}k^{p-1}\Bigr]_+ -3L\rho_0>\eps.$$
Then, choosing 
$$\eta=c_1\bigl(\rho^2\tau^{-1}+1\bigr)^{-1/(1-\beta)},$$
with $c_1=c_1(p,n,k)\in (0,\eps)$ sufficiently small, we get
$$\begin{aligned}
(1-\beta)&(1-C_2\eta)_+^2\Bigl[\beta-(1+\eps)(1+C_2\eta)^{p-2}k^{p-1}\Bigr]_+
-3L\rho-C_5\eta^{1-\beta}\bigl(\rho^2\tau^{-1}+1+\eps^{1-p}\bigr)\\
&\ge \eps-C_5\eta^{1-\beta}\bigl(\rho^2\tau^{-1}+1+\eps^{1-p}\bigr)\ge 0.
\end{aligned}$$
 Hence,
$$V_t-\Delta V-|\nabla V|^p\ge 0
\quad\hbox{ in $D\times (t_0,t_1)$.}$$
On the other hand, we have
$v=0=V$ on $(B_\rho(a)\cap\partial\omega)\times (t_0,t_1)$
and, by assumption \eqref{HypTypeILemSmall1},
$$v\le \kappa\delta^{1-\beta}=V \quad\hbox{ in $\bigl[D\times \{t_0\}\bigr]\cup\bigl[(\omega\cap\partial B_\rho(a))\times (t_0,t_1)\bigr]$}.$$
It then follows from the comparison principle that $v\le V$ in $D\times (t_0,t_1)$. In particular,
for all $(x,t)\in(B_{\rho/2}(a)\cap\partial\omega)\times (t_0,t_1)$, we obtain
$$\begin{aligned}
v_\nu\le V_\nu
&=(1-\beta)\kappa\varphi^{-\beta}=k(\eta\rho)^{-\beta} \bigl[\tau^{-1}(t-t_0)\bigr]^{-\frac{\beta}{1-\beta}} \\
&=kc_1^{-\beta}\rho^{-\beta}\bigl(1+\rho^2\tau^{-1}\bigr)^{\frac{\beta}{1-\beta}} \tau^{\frac{\beta}{1-\beta}}(t-t_0)^{-\frac{\beta}{1-\beta}}
=kc_1^{-\beta}\rho^{-\beta}\left(\frac{\tau+\rho^2}{t-t_0}\right)^{\frac{\beta}{1-\beta}},
\end{aligned}$$
which is the desired conclusion.
\qed

\medskip

{\it Proof of Lemma \ref{lem:exclude blow-up points}.}
Put $w=|\nabla v|^2$, then
$$
\mathcal{L}w=-2|D^2v|^2,
$$
where $|D^2v|^2=\Sigma_{i,j}(v_{x_ix_j})^2$ and
$$
\mathcal{L}w=w_t-\Delta w-p|\nabla v|^{p-2}\nabla v\cdot\nabla w.
$$
Let $m\in(0,1)$. We select a cut-off function
$\eta\in C^2(\overline{B}_{R}(x_0))$, with
$$\eta=0\ \ {\rm for}\ |x-x_0|=R,\quad 0<\eta\leq 1\ \ {\rm in}\ B_{R}(x_0)$$
and such that
\begin{equation*}
    |\nabla\eta|\leq CR^{-1}\eta^m,\ \ \
    |\Delta\eta|+4\eta^{-1}|\nabla\eta|^2\leq CR^{-2}\eta^m
    \ \ {\rm in}\ B_{R}(x_0),
\end{equation*}
 where $C=C(m)>0$. Such a function $\eta$ is given for instance in the proof of
 \cite[Theorem~3.2]{SZ}. 
 Put
$$z=\eta w,\quad (x,t)\in Q:=(B_{R}(x_0)\cap\omega)\times(t_0,t_1).$$
Then
\begin{equation*}
    \mathcal{L}z=\eta\mathcal{L}w+w\mathcal{L}\eta-2\nabla\eta\cdot\nabla w  \ \ {\rm in}\ Q.
\end{equation*}
Since $2|\nabla\eta\cdot\nabla w|\leq 4\eta^{-1}|\nabla
\eta|^2w+\eta|D^2u|^2$, it follows that
$$    \mathcal{L}z+\eta|D^2u|^2\leq w(|\Delta\eta|+4\eta^{-1}|\nabla\eta|^2)+pw^{(p+1)/2}|\nabla\eta|
                      \leq CR^{-2}\eta^m w+CpR^{-1}\eta^m w^{(p+1)/2}.
$$
Using $|w^{p/2}-u_t|=|\Delta u|\leq \sqrt{n}|D^2u|^2$,
hence $w^p/(2n)\leq |D^2u|^2+|u_t|^2$, we get
\begin{equation*}
    \mathcal{L}z+\frac{1}{2n}\eta w^p\leq CR^{-2}\eta^m
    w+CpR^{-1}\eta^m w^{(p+1)/2}+M^2.
\end{equation*}
Taking $m=(p+1)/(2p)$ and using Young's inequality, we obtain
$$
    \mathcal{L}z+\frac{1}{4n}z^p\leq
     \mathcal{L}z+\frac{1}{4n}\eta w^p\leq A:=CR^{-2p/(p-1)}+ M^2
    \ \ \ {\rm in}\ \ (B_R(x_0)\cap\omega)\times(t_0,t_1).
$$
Let $c=[4n/(p-1)]^{1/(p-1)}$ and set
$$\phi:=B+c(t-t_0)^{-1/(p-1)}\ \hbox{ for $t>t_0$, \quad where $B=\max\bigl\{N^2,(4nA)^{1/p}\bigr\}$.}$$
We have $\phi\ge B\ge z$ on $\partial((B_R(x_0)\cap\omega)\times(t_0,t_1)$, and
$$\begin{aligned}
\mathcal{L}\phi+\frac{1}{4n}\phi^p
&=-\frac{c}{p-1}(t-t_0)^{-p/(p-1)}+\frac{1}{4n} \bigl(B+c(t-t_0)^{-1/(p-1)}\bigr)^p \\
&\ge \Bigl[-\frac{c}{p-1}+\frac{c^p}{4n}\Bigr](t-t_0)^{-p/(p-1)}+\frac{B^p}{4n}\ge A
 \ \ \ {\rm in}\ \ (B_R(x_0)\cap\omega)\times(t_0,t_1).
\end{aligned}$$
Since $z$ remains bounded as $t\to t_0$, whereas $\phi(t)\to \infty$,
it follows from the comparison principle (cf.~e.g.~\cite[Proposition 2.2]{SZ}) that
$\phi\ge z$ in $(B_R(x_0)\cap\omega)\times(t_0,t_1)$. Hence in particular
$$\begin{aligned}
    |\nabla v|&\leq B^{1/2}+c^{1/2}(t-t_0)^{-1/2(p-1)}   \ \ \ {\rm in}\ \ (B_{R/2}(x_0)\cap\omega)\times(t_0,t_1),
\end{aligned}$$
which implies the desired conclusion.
\qed

\medskip

{\it Proof of Proposition \ref{lemCentral}.} Assumption \eqref{innersphere} guarantees that,
for each $b\in B_\sigma(a)\cap\partial\omega$, the line segment $(b,b+R\nu_b)$ lies in $D$
and that
\be\label{innersphere2}
{\rm dist}(b+\rho\nu_b,\partial\omega)=\rho,\quad 0<\rho<R.
\ee
Let $r_0=r_0(p,k,M,C_\eps)$ be given by Lemma \ref{lem2b}. We deduce from assumptions
\eqref{boundut1central}-\eqref{HypTypeILemSmall1central} and Lemma \ref{lem2b} that
$$|v_\nu(b+\rho\nu_b,t)| \le k\rho^{-\beta},
\quad\hbox{ for all $b\in B_\sigma(a)\cap\partial\omega$, all $\rho\in (0,r]$ and all $t\in (T_0-\theta,T_0)$.}$$
Since $v=0$ on $ B_\sigma(a)\cap\partial\omega$, it follows by integration that
$$|v(b+\rho\nu_b,t)| \le (1-\beta)^{-1}k\rho^{1-\beta} 
\quad\hbox{ for all $b\in B_\sigma(a)\cap\partial\omega$, all $\rho\in (0,r]$ and all $t\in (T_0-\theta,T_0)$.}$$
Using \eqref{innersphere2} and $\sigma\le r\le R$, we have in particular
$$|v(x,t)| \le (1-\beta)^{-1}k\,{\rm dist}^{1-\beta}(x,\partial\omega)  
\quad  \hbox{말n $(\omega\cap B_{\sigma/2}(a))\times (T_0-\theta,T_0)$}.$$
We may then apply Lemma \ref{lem0bApp} (taking $r_0$ possibly smaller, which may also depend on $L$) 
to infer that
$$
\partial_\nu v \le c(n,p,k)\sigma^{-\frac{1}{p-1}}\bigl[1+\sigma^2\theta^{-1}\bigr]^{\frac{1}{p-2}} 
\quad\hbox{ in $(B_{\sigma/4}(a)\cap\partial\omega)\times (T_0-\theta/2,T_0)$}.
$$
Since $z:=-v$ satisfies $z_t-\Delta z=-|\nabla z|^p\le |\nabla z|^p$, we may apply Lemma \ref{lem0bApp} 
to $z$ as well, so that actually
$$
|\partial_\nu v|\le c(n,p,k)\sigma^{-\frac{1}{p-1}}\bigl[1+\sigma^2\theta^{-1}\bigr]^{\frac{1}{p-2}}
\quad\hbox{ in $(B_{\sigma/4}(a)\cap\partial\omega)\times (T_0-\theta/2,T_0)$}.
$$
From this estimate, \eqref{ConclLemCentral} finally follows from Lemma \ref{lem:exclude blow-up points}.
\qed

\bigskip

{\bf Acknowledgements.} Part of this work was done during a visit of
PhS at the
Dipartimento di Matematica e Informatica of the Universit\`a degli Studi di Perugia
 within the auspices of the INdAM -- GNAMPA Projects
2018. He wishes to thank this institution for the kind hospitality.
PhS is partly supported by the Labex MME-DII (ANR11-LBX-0023-01). 
\smallskip

RF and PP were partly supported by the Italian MIUR project
{\em Variational methods, with applications to problems in mathematical physics and 
geometry (2015KB9WPT\_009) and are members of the {\em Gruppo Nazionale per
l'Analisi Ma\-te\-ma\-ti\-ca, la Probabilit\`a e le loro Applicazioni}
(GNAMPA) of the {\em Istituto Nazionale di Alta Matematica} (INdAM).
The manuscript was realized within the auspices of the INdAM -- GNAMPA Projects
2018
{\em Pro\-ble\-mi non lineari alle derivate parziali} (Prot\_U-UFMBAZ-2018-000384).}
\medskip

{\bf Declaration of interest statement.} No conflict of interests.


\begin{thebibliography}{99}
\bibitem{Alaa}
{\sc N. Alaa}, 
Weak solutions of quasilinear parabolic equations with measures as initial data,
Ann. Math. Blaise Pascal 3 (1996) 1--15.

\bibitem{AGQ}
{\sc S. Alarc\'on, J. Garc\'\i a-Meli\'an, A. Quaas},
Keller-Osserman type conditions for some elliptic problems with gradient terms,
J. Differential Equations 252 (2012) 886--914.

\bibitem{ABG}
{\sc N.D. Alikakos, P.W. Bates, C.P. Grant},
Blow up for a diffusion-advection equation,
 Proc. Roy. Soc. Edinburgh Sect. A 113 (1989) 181--190.

\bibitem{Ar69} {\sc D.G. Aronson},
Regularity properties of flows through porous media,
SIAM J. Applied Math. 17 (1969) 461--467.

\bibitem{ARS04}
{\sc J.M. Arrieta, A. Rodr\'\i guez-Bernal, Ph. Souplet},
Boundedness of global solutions for nonlinear parabolic equations
involving gradient blow-up phenomena,
Ann. Scuola Norm. Sup. Pisa Cl. Sci. (5)
3 (2004) 1--15.

\bibitem{BG}
{\sc C. Bandle, E. Giarrusso},
Boundary blow-up for semilinear elliptic equations with nonlinear gradient terms,
 Adv. Differential Equations 1 (1996) 133--150.

\bibitem{BaLio04}
{\sc G. Barles, F. Da Lio},
On the generalized Dirichlet problem for viscous Hamilton-Jacobi equations,
J. Math. Pures Appl. 83 (2004) 53--75.

\bibitem{BCN}
{\sc H. Berestycki, L. Caffarelli, L. Nirenberg},
Further qualitative properties for elliptic equations in unbouded domains,
Ann. Sc. Norm. Super. Pisa Cl. Sci.
(4) 15 (1997) 69--94.

\bibitem{Ber}
{\sc S. Bernstein},
Sur la g\'en\'eralisation du probl\`eme de Dirichlet,
Math. Ann. 69 (1910) 82--136.

\bibitem{CG96}
{\sc G. Conner, C. Grant},
Asymptotics of blowup for a convection-diffusion equation with conservation,
Differential Integral Equations 9 (1996) 719--728.

\bibitem{CrM}
{\sc G. Crasta, A. Malusa}, 
The distance function from the boundary in a Minkowski space,
Trans. Amer. Math. Soc. 359 (2007), 5725--5759.

\bibitem{Dup}
{\sc L. Dupaigne},
Stable solutions of elliptic partial differential equations.
Monographs and Surveys in Pure and Applied Mathematics, 143,
Chapman $\&$ Hall/CRC, Boca Raton, FL, 2011, xiv+321 pp.

\bibitem{Esteve19}
{\sc C. Esteve},
Single-point gradient blow-up on the boundary for diffusive Ham\-ilton-Jacobi equation in domains with non-constant curvature,
Pre\-print arXiv:1902.03080
(2019).

\bibitem{Fa15}
{\sc A. Farina},
Some symmetry results and Liouville-type theorems for solutions to semilinear equations,
Nonlinear Anal. 121 (2015) 223--229.

\bibitem{FV}
{\sc A. Farina, E. Valdinoci},
Flattening results for elliptic PDEs in unbounded domains with applications to overdetermined problems,
Arch. Ration. Mech. Anal. 195 (2010) 1025--1058.

\bibitem{FQS}
{\sc P. Felmer, A. Quaas, B. Sirakov},
Solvability of nonlinear elliptic equations with gradient terms,
J. Differential Equations 254 (2013) 4327--4346.

\bibitem{fl19}
{\sc M. Fila, J. Lankeit},
Continuation beyond interior gradient blow-up in a semilinear parabolic equation,
Preprint arXiv:1902.01127 (2019).

\bibitem{FL94}
{\sc M. Fila, G. Lieberman},
Derivative blow-up and beyond for quasilinear parabolic equations,
 Differential Integral Equations 7 (1994) 811--821.

\bibitem{fi_aihp14}
{\sc Y. Fujishima, K. Ishige},
Blow-up set for type I blowing up solutions for a semilinear heat equation,
Ann. Inst. H. Poincar\'e Anal. Non Lin\'eaire
31 (2014) 231--247.

\bibitem{GK87}
{\sc Y. Giga, R.V. Kohn},
Characterizing blowup using similarity variables,
Indiana Univ. Math. J.
36 (1987) 1--40.

\bibitem{GuoHu}
{\sc J.-S. Guo, B. Hu},
Blowup rate estimates for the heat equation with a nonlinear gradient source term
Discrete Contin. Dyn. Syst. 20 (2008) 927--937.

\bibitem{GuoSou2}
{\sc J.-S. Guo, Ph. Souplet},
Excluding blowup at zero points of the potential by means of Liouville-type theorems,
J. Differential Equations 265 (2018) 4942--4964.

\bibitem{HM} 
{\sc M. Hesaaraki, A. Moameni},
Blow-up positive solutions for a family of nonlinear parabolic equations in general domain in $R^N$,
Mich. Math. J. 52 (2004) 375--389.

\bibitem{La58} {\sc O.A. Lady\v{z}enskaja},
Solution of the first boundary problem in the large for quasi-linear parabolic equations,
Trudy Moskov. Mat. Ob\v s\v c. 7 (1958) 149--177.

\bibitem{LL}
{\sc J.M. Lasry, P.L. Lions},
Nonlinear elliptic equations with singular boundary conditions and stochastic control with state constraints.
I. The model problem, Math. Ann. 283 (1989) 583--630.

\bibitem{LP1}
{\sc T. Leonori, A. Porretta},
The boundary behavior of blow-up solutions related to a stochastic control problem with state constraint,
SIAM J. Math. Anal. 39 (2007/08) 1295--1327.

\bibitem{LP2}
{\sc T. Leonori, A. Porretta},
Gradient bounds for elliptic problems singular at the boundary,
Arch. Ration. Mech. Anal. 202 (2011) 663--705.

\bibitem{Lions85}
{\sc P.L. Lions},
Quelques remarques sur les probl\`emes elliptiques quasilin\'eaires du second ordre,
J. Anal. Math.  45 (1985) 234--254.

\bibitem{LY75}
P. Li, S.-T. Yau,
On the parabolic kernel of the Schr\"odinger operator,
Acta Math. 156 (1986) 153--201.

\bibitem{LS10}
{\sc Y.-X. Li, Ph. Souplet},
Single-point gradient blow-up on the boundary for diffusive Hamilton-Jacobi equations in planar domains,
Comm. Math. Phys. 293 (2010) 499--517.

\bibitem{MZ98}
{\sc F. Merle, H. Zaag},
Optimal estimates for blowup rate and behavior for nonlinear heat equations,
Comm. Pure Appl. Math. 51 (1998) 139--196.

\bibitem{MZ00}
{\sc F. Merle, H. Zaag},
A Liouville theorem for vector-valued nonlinear heat equations and applications,
Math. Ann. 316 (2000) 103--137.

\bibitem{MizNeu}
{\sc N. Mizoguchi},
Blowup rate of solutions for a semilinear heat equation with the Neumann boundary condition,
J. Differential Equations
193 (2003) 212--238.

\bibitem{PS_imrn16}
{\sc A. Porretta, Ph. Souplet},
The profile of boundary gradient blowup for the diffusive Hamilton-Jacobi equation,
Int. Math. Res. Notices 17 (2017) 5260--5301.

\bibitem{PS_aihp17}
{\sc A. Porretta, Ph. Souplet},
Analysis of the loss of boundary conditions for the diffusive Hamilton-Jacobi equation,
Ann. Inst. H. Poincar\'e Anal. Non Lin\'eaire 34 (2017) 1913--1923.

\bibitem{PS_jmpa18}
{\sc A. Porretta, Ph. Souplet},
Blow-up and regularization rates, loss and recovery of boundary conditions
 for the superquadratic viscous Hamilton-Jacobi equation,
J. Math. Pures Appl. doi.org/10.1016/j.matpur.2019.02.014,
to appear (Preprint arXiv:\allowbreak1811.01612).

\bibitem{PV}
{\sc A. Porretta, L. V\'eron},
Asymptotic behavior for the gradient of large solutions to some nonlinear elliptic equations,
Adv. Nonlinear Studies   6 (2006)  351--378.

\bibitem{PZ_aihp12}
{\sc A. Porretta, E. Zuazua},
Null controllability of viscous Hamilton-Jacobi equations
Ann. Inst. H. Poincar\'e Anal. Non Lin\'eaire 29 (2012) 301--333.

\bibitem{QR18}
{\sc A. Quaas, A. Rodr\'\i guez},
Loss of boundary conditions for fully nonlinear parabolic equations with superquadratic gradient terms,
J. Differential Equations 264 (2018) 2897--2935.

\bibitem{Se69} {\sc J. Serrin}, The problem of Dirichlet for quasilinear elliptic
differential equations with many
independent variables, Phil. Trans. Roy. Soc. London A 264 (1969) 413--469.

\bibitem{Sou} 
{\sc Ph. Souplet},
Gradient blow-up for multidimensional nonlinear parabolic equations with general boundary conditions,
 Differential Integral Equations 15 (2002) 237--256.

\bibitem{SZ}
{\sc P. Souplet, Q.S. Zhang},
Global solution of inhomogeneous Hamilton-Jacobi equations, J. Anal. Math. 99 (2006) 355--396.

\bibitem{ZL13}
{\sc Z.-C. Zhang, Z. Li},
A note on gradient blowup rate of the inhomogeneous Hamilton-Jacobi equations,
Acta Math. Sci. Ser. B (Engl. Ed.)
33 (2013) 678--686.
\end{thebibliography}
\end{document}